\documentclass[11pt]{amsart}
\usepackage{geometry}                
\usepackage[parfill]{parskip}    
\usepackage{graphicx}
\usepackage{amssymb}
\usepackage{epstopdf}
\usepackage{subfig}
\usepackage{amsmath}
\usepackage{color}
\usepackage{pdfsync}

\swapnumbers 
\numberwithin{equation}{section} 

\theoremstyle{plain} 
\newtheorem{thm}{Theorem}[section]

\newtheorem{lem}[thm]{Lemma} 
\newtheorem{conjec}[thm]{Conjecture} 
\newtheorem{prop}[thm]{Proposition} 
\allowdisplaybreaks

\graphicspath{/4-3-2012-Moscow}
\theoremstyle{definition} 
\newtheorem{defin}[thm]{Definition}

\newtheorem{rem}[thm]{Remark}

\def\pp{\medskip{\parindent 0pt \it Proof.\ }} 
\graphicspath{/4-3-2012-Moscow/}

\title[Vassiliev Invariants of Virtual Legendrian Knots]{Vassilev Invariants of Virtual Legendrian Knots}
\author[Cahn \& Levi]{Patricia Cahn and Asa Levi} \address{Patricia Cahn, Department of Mathematics, Department of Mathematics, David Rittenhouse Lab.
209 South 33rd Street, Philadelphia, PA 19104-6395, USA} \email{pcahn@math.upenn.edu} \address{Asa Levi, Department of Mathematics, 6188 Kemeny Hall, Dartmouth College, Hanover NH 03755-3551, USA} \email{Asa.Levi.gr@dartmouth.edu} 
\date{}

\begin{document}
\maketitle
\begin{abstract} We introduce a theory of virtual Legendrian knots.  A virtual Legendrian knot is a cooriented wavefront on an oriented surface up to Legendrian isotopy of its lift to the unit cotangent bundle and stabilization and destablization of the surface away from the wavefront. We show that the groups of Vassiliev invariants of virtual Legendrian knots and of virtual framed knots are isomorphic.  In particular, Vassiliev invariants cannot be used to distinguish virtual Legendrian knots that are isotopic as virtual framed knots and have equal virtual Maslov numbers.
\end{abstract}
We work in the smooth category.  All maps and manifolds are $C^\infty$. All surfaces are oriented unless explicitly stated otherwise.  All curves are immersed.

\section{Introduction}

  The first goal of this paper is to introduce a theory of virtual Legendrian knots. Briefly, a virtual Legendrian knot is a Legendrian knot in the spherical cotangent bundle $ST^*F$ of a surface $F$, up to Legendrian isotopy and stabilization and destabilization of the surface $F$.  In this paper, virtual framed knots are framed knots in $ST^*F$ up to framed isotopy and stabilization and destabilization of the surface $F$. The second goal of the paper is to study the group of order $\leq n$ Vassiliev invariants of virtual Legendrian knots.  We show that this group is naturally isomorphic to the group of order $\leq n$ Vassiliev invariants of virtual framed knots.  As a corollary, order $\leq n$ Vassiliev invariants do not distinguish virtual Legendrian knots that are isotopic as framed virtual knots and have the same virtual Maslov number.  Similar theorems were proved by Goryunov \cite{Goryunov}, Hill \cite{Hill}, Fuchs-Tabachnikov \cite{f&t}, and Chernov \cite{Chernov} for Legendrian knots in various contact manifolds.

We give three equivalent formulations of virtual Legendrian knot theory.  The first, which we describe in the introduction, was motivated by the formulation of virtual knot theory given by Carter, Kamada and Saito \cite{CKS}, and was suggested to us by Chernov \cite{ChernovDefn}.

First we review the concept of a wavefront as discussed by Arnold in \cite{ArnoldMechanics}. Suppose that the surface $F$ is made of an isotropic, homogeneous medium, and that light rays are emitted from a point of $F$.  The set of points that these light rays reach at a fixed time $t$ is called a wavefront.  As $t$ increases, semi-cubical cusps may appear, so the wavefront is not necessarily immersed.  We say a (cooriented) {\it wavefront} on an oriented surface $F$ is a cooriented curve on $F$ which is immersed except at a finite number of semi-cubical cusps.  The coorientation represents the direction of propagation of the wavefront.  A wavefront is {\it generic} if it has a finite number of self-intersection points, all of which are transverse double points.

\begin{figure}[htbp]

	\includegraphics[width=5cm]{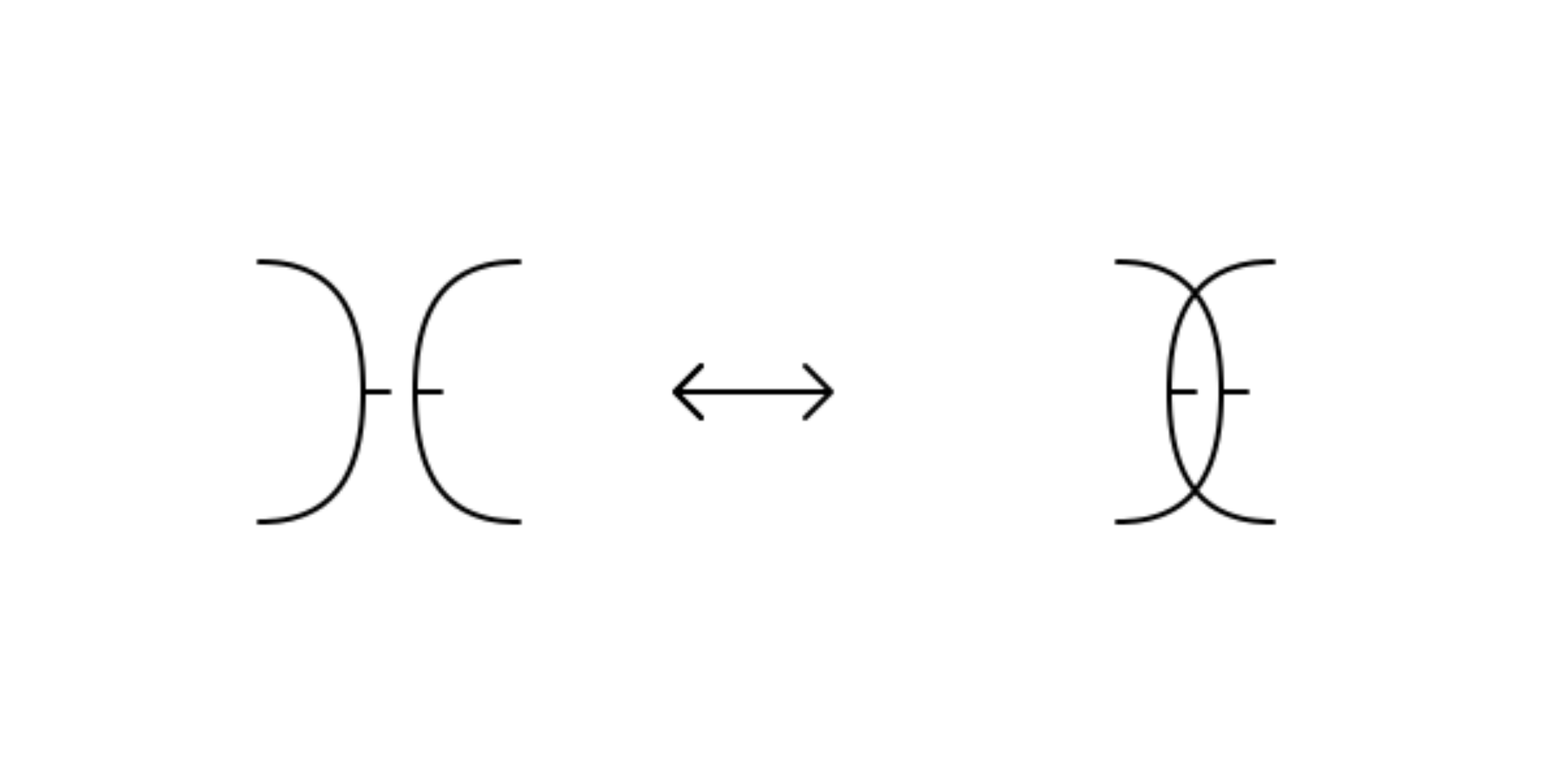}
	\caption{The dangerous tangency move.}
	\label{dangerous.fig}
\end{figure}

The spherical cotangent bundle $ST^*F$ of a surface $F$ is equipped with the natural contact structure.  Arnold observed that, due to the Huygens principle, the propagation of a wavefront on $F$ lifts to a Legendrian isotopy in $ST^*F$ \cite{ArnoldMechanics, ArnoldPDE}.  That is, lifting the wavefront to $ST^*F$ according to the direction of its coorienting vector produces a curve in $ST^*F$ which is everywhere tangent to the distribution of contact planes, and during the propagation of the wavefront, the lift of this curve undergoes an isotopy while remaining tangent to the contact planes.  

In particular the Huygens principle implies that the {\it dangerous tangency move} in Figure \ref{dangerous.fig} cannot appear during the propagation of a single front $\omega$, because if two branches of $\omega$ become tangent during the propagation in such a way that their coorientations match, they must be tangent for all $t$.

\begin{figure}[htbp]
	\includegraphics[width=6cm]{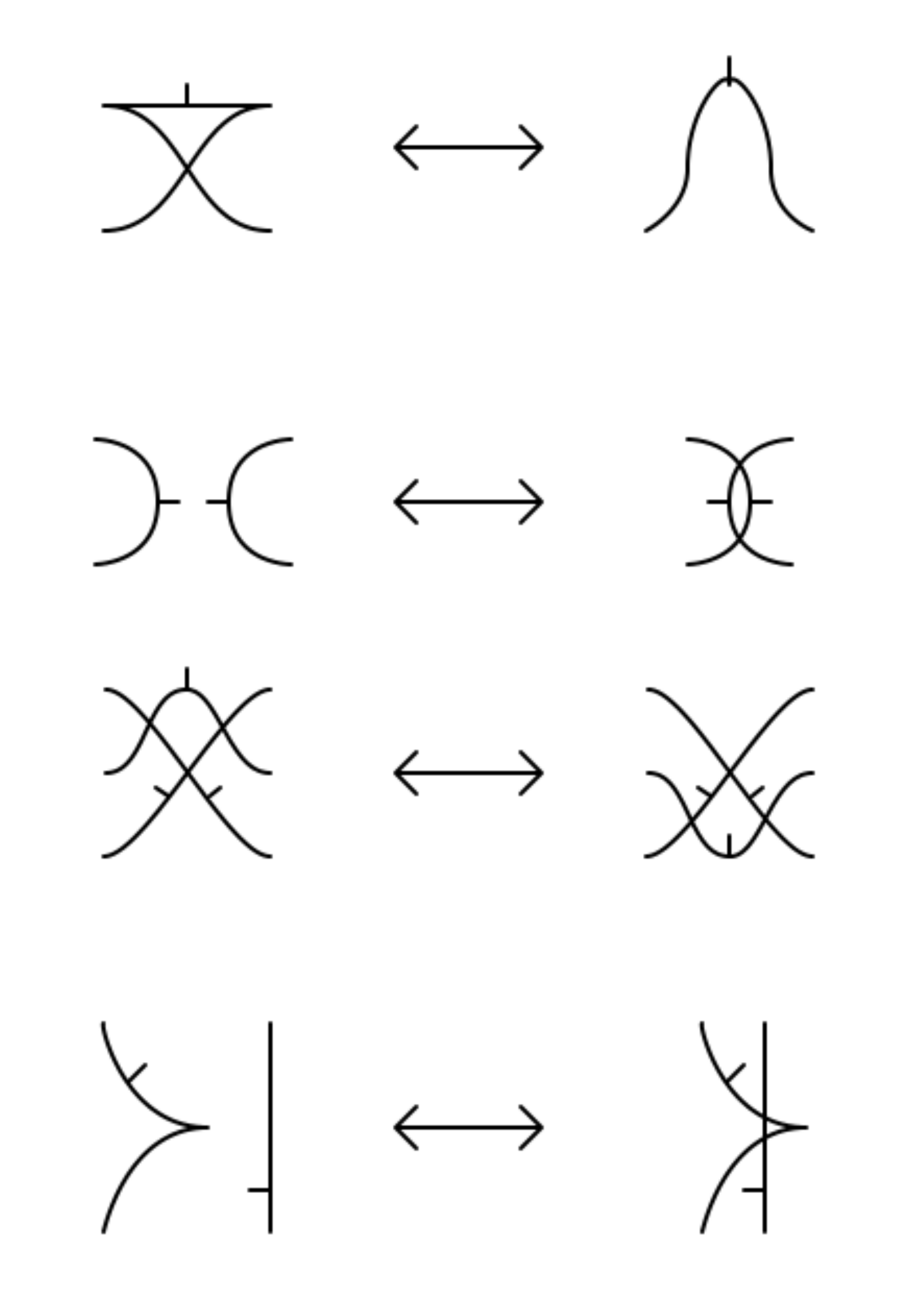}
	\caption{Moves for wavefronts on an oriented surface, with one possible choice of coorientation on each branch.}
	\label{wavefrontmoves.fig}
\end{figure}

Hence the Legendrian liftings of two generic wavefronts are Legendrian isotopic if and only if the wavefronts are related by a sequence of the moves in Figure \ref{wavefrontmoves.fig} up to certain choices of coorientation, in addition to ambient isotopy.  To obtain all the valid choices of coorientation, one should consider the moves in Figure \ref{wavefrontmoves.fig} with all possible choices of coorientations on the branches, {\it except} in the case of the second move, because dangerous tangencies are prohibited \cite{ArnoldInvariants}.

A virtual Legendrian knot is a Legendrian knot in the spherical cotangent bundle of a surface, and corresponds to a wavefront on that surface.  The definition of virtual Legendrian knot theory suggested by Chernov is as follows: two virtual Legendrian knots are equivalent if their corresponding wavefronts are related by wavefront moves, and stabilization and destablization of the surface. To stabilize the surface $F$, we remove two disks from $F$ that are disjoint from the wavefront, and glue the two boundary components together so that the surface remains orientable.  To destabilize the surface $F$, we choose an essential simple closed curve disjoint from the wavefront, cut along it, and glue disks to the two resulting boundary components.  Physically, this corresponds to the notion that the medium through which the wavefront propagates can change its topology outside a neighborhood of the wavefront.  We denote a virtual Legendrian knot by a pair $(F,K)$ where $F$ is a compact oriented surface and $K$ is a wavefront on $F$, and we let $[(F,K)]_l$ denote its virtual Legendrian isotopy class.

We also give a purely diagrammatic description of virtual Legendrian knot theory in Section \ref{diagrams}.  Namely, one can view a virtual Legendrian knot as a wavefront in the plane with virtual crossings up to certain moves, in the spirit of Kauffman's orginal theory of virtual knots \cite{Kauffman}.

Questions similar to ours have been studied by many people over the last 20 years.  In the rest of this section let $\mathcal{K}$ be $\mathbb{R}$ or $\mathbb{C}$ and $\mathcal{A}$ be any abelian group.  In 1997 Fuchs and Tabachnikov \cite{f&t} showed that the vector spaces of $\mathcal{K}$-valued, order $\leq n$ Vassiliev invariants of Legendrian knots with fixed Maslov number in $\mathbb{R}^3$ endowed with the standard contact structure and of framed knots in $\mathbb{R}^3$ are isomorphic.  Around the same time, Goryunov showed that the vector spaces of $\mathcal{K}$-valued, order $\leq n$ Vassiliev invariants of oriented framed knots in the solid torus $ST^*\mathbb{R}^2$ and of oriented plane curves without direct self-tangencies are isomorphic \cite{Goryunov}.

Generalizing the work of Goryunov, Hill \cite{Hill} proved the same result for all planar fronts.  Finally, Chernov \cite{Chernov} was able to show that the groups of $\mathcal{A}$-valued, order $\leq n$ Vassiliev invariants of framed knots in $ST^*F$ where $F$ is any surface, and of Legendrian knots in $ST^*F$ are isomorphic.

More precisely, suppose $\mathcal{L}$ is a connected component of the space of Legendrian curves in a contact manifold $(M,C)$ and let $\mathcal{F}$ be the connected component of the space of framed curves in $(M,C)$ that contains $\mathcal{L}$.  Chernov proved that the groups of $\mathcal{A}$-valued, order $\leq n$ Vassiliev invariants of framed knots in $\mathcal{F}$ and Legendrian knots in $\mathcal{L}$ are isomorphic for a large class of contact manifolds $(M,C)$.  In particular, the theorem holds for $M=ST^*F$ with its natural contact structure. We develop techniques to show that a similar result holds in the virtual category.
\begin{thm} \label{mainthm}
Let $\mathcal{F}$ be a connected component of the space of virtual framed curves and $\mathcal{L}\subset\mathcal{F}$ be a connected component of the space of virtual Legendrian curves contained in $\mathcal{F}$.  Let $\mathcal{A}$ be an abelian group and $\mathcal{V}_n^\mathcal{F}$ (respectively $\mathcal{V}_n^\mathcal{L}$) be the group of $\mathcal{A}$-valued Vassiliev invariants of framed (respectively Legendrian) knots on $\mathcal{F}$ (respectively $\mathcal{L}$) of order $\leq n$.   Then the restriction map $\phi:\mathcal{V}_n^{F}\rightarrow\mathcal{V}_n^\mathcal{L}$ is an isomorphism.
\end{thm}
It follows that Vassiliev invariants cannot distinguish two virtual Legendrian knots that are homotopic as Legendrian curves and isotopic as framed virtual knots.

\begin{thm}\label{maincor}
Let $x\in \mathcal{V}_n^\mathcal{L}$, let $(F, K), (F', L)$ be representatives of the virtual Legendrian homotopy class $\mathcal{L}$, and let $(F, K)$ be virtually framed isotopic to $(F', L')$.  Then $x([(F, K)]_l) = x([(F', L)]_l)$.
\end{thm}

In Section \ref{classical.sec} we discuss virtual version of the Maslov number.  We prove that, as in the classical case, virtual Legendrian homotopy classes are completely characterized by their Maslov number and their underlying virtual homotopy class.

\begin{thm} Two virtual Legendrian knots are virtual Legendrian homotopic if and only if they have the same virtual Maslov number and are homotopic as virtual knots.
\end{thm}

Part of the motivation for studying virtual knot theory stems from the fact that virtually isotopic classical knots must be isotopic as classical knots.  This was first proven by Goussarov, Polyak, and Viro \cite{GPV}, and the proof can also be found in \cite{Kauffman}.  In other words, virtual knot theory extends classical knot theory.

Together with Chernov, we conjecture:
\begin{conjec}  Two Legendrian knots in $ST^*F$, with $F=S^2$ or $\mathbb{R}^2$, that are isotopic as virtual Legendrian knots must be Legendrian isotopic in $ST^*F$.
\end{conjec}
  
Kuperberg \cite{Kuperberg} later showed that two knots in $F\times I$ that are isotopic as virtual knots must be isotopic as knots in $F\times I$, possibly after a homeomorphism of $F\times I$, provided that $F$ is the surface of smallest genus realizing an element of the virtual isotopy class.  We also hope that a similar result holds for virtual Legendrian knots.

\begin{conjec} Let $K_1$ and $K_2$ be two Legendrian knots in $ST^*F$ that are isotopic as virtual Legendrian knots, and suppose that $F$ is the surface of smallest genus realizing knots in the virtual Legendrian isotopy class of $K_1$ and $K_2$.  Then, possibly after a contactomorphism of $ST^*F$, $K_1$ and $K_2$ are Legendrian isotopic in $ST^*F$.
\end{conjec}

Chernov's definition of virtual Legendrian knot theory can also be generalized to higher dimensions.  That is, two Legendrian manifolds $L_1$ in $ST^*M_1$ and $L_2$ in $ST^*M_2$ (not necessarily spheres) are virtually Legendrian isotopic if one can be obtained from the other by a sequence of Legendrian isotopies and modifications of contact $ST^*M $ induced by surgery on $M$ in the part of $M$ which does not contain the projection of the Legendrian knot.  It's also possible to formulate the first conjecture in higher dimensions, namely two Legendrian knots in $ST^*\mathbb{R}^m$ or $ST^*S^m$ are virtual Legendrian isotopic if and only if they are Legendrian isotopic.

\section{Virtual Legendrian Knot Diagrams}\label{diagrams}
A {\it virtual Legendrian knot diagram} is a generic wavefront in $\mathbb{R}^2$ with two types of crossings.  In addition to ordinary crossings, there are virtual crossings, which are marked with a small circle and obey slightly different Reidemeister moves.  The Reidemeister moves involving only ordinary crossings are shown in Figure \ref{wavefrontmoves.fig}, and these moves are a subset of the possible moves for virtual Legendrian knot diagrams.  Again, we can obtain other wavefront moves from these moves by independently reversing the choice of coorientation on any branch, except in the case of the second Reidemeister move, where the dangerous tangency move is forbidden.

\begin{figure}[htbp]
		\includegraphics[width=8cm]{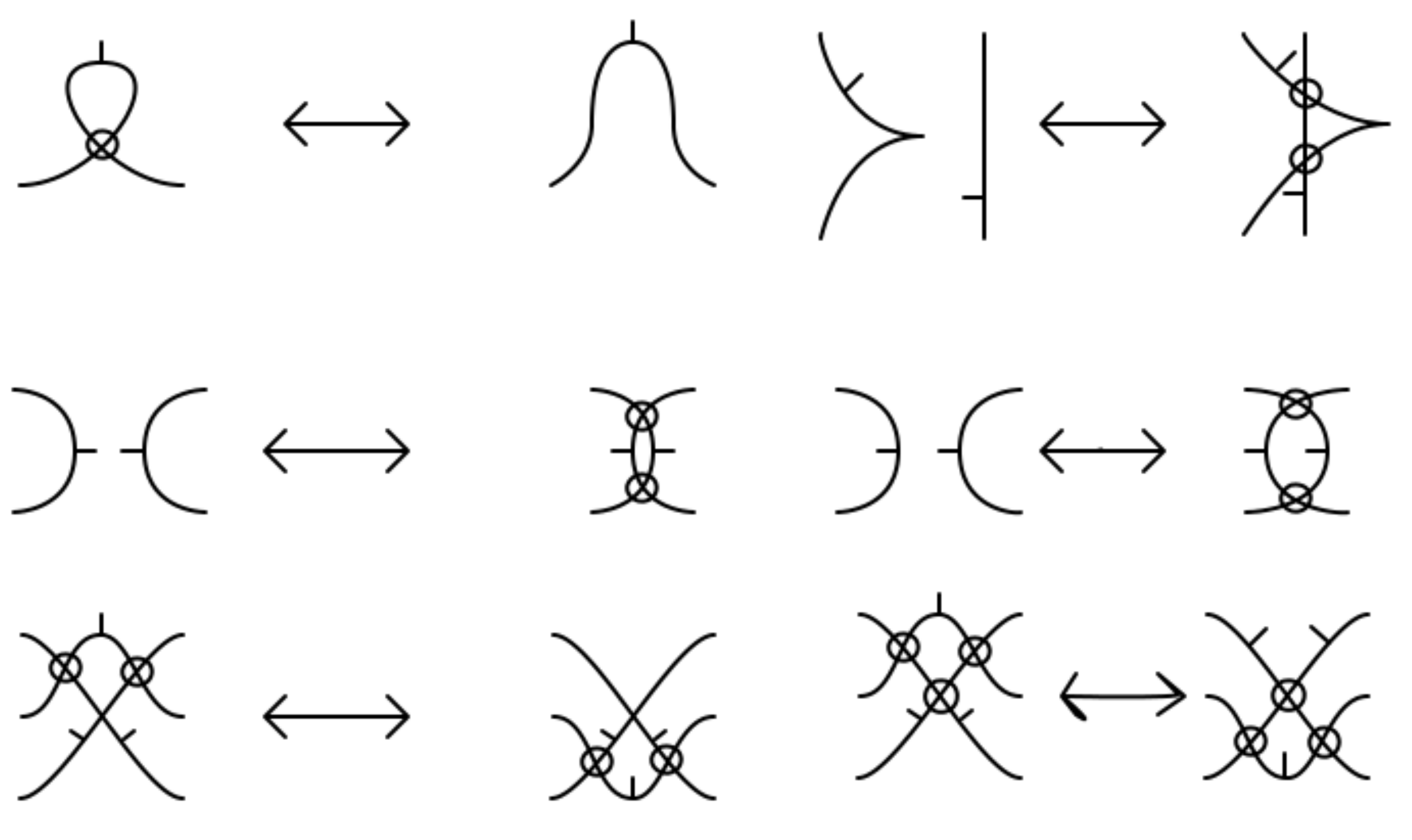}
		\caption{Moves for virtual wavefront diagrams in the plane, with one possible choice of coorientation on each branch.}
		\label{virtualwavefrontmoves.fig}
\end{figure}

The remaining moves for virtual Legendrian front diagrams involve at least one virtual crossing, and are pictured in Figure \ref{virtualwavefrontmoves.fig}.  Other moves can be obtained from these moves by independently reversing the coorientation on any branch.  Note that we allow virtual dangerous tangencies. Also note that the move obtained from the first move in Figure \ref{wavefrontmoves.fig} by replacing the ordinary crossing with a virtual crossing is {\it not} allowed.  

In Section \ref{equivalence} we will verify that this diagrammatic definition of virtual Legendrian knot theory is equivalent to Chernov's definition given in the Introduction.

If we allow the dangerous self-tangency move in Figure \ref{dangerous.fig}, the equivalence relation generated is that of {\it virtual Legendrian homotopy}, rather than isotopy.  We sometimes refer to a virtual Legendrian homotopy class as a connected component of the space of virtual Legendrian curves.

\section{Flat Virtual Knots}\label{virtualstrings}
Virtual knot theory was introduced by Kauffman \cite{Kauffman}.  We will briefly review the definition of a related object, called a {\it flat virtual knot}, or {\it virtual string}, in order to motivate the definition of a virtual Legendrian knot.  Virtual strings were introduced by Turaev \cite{Turaev}.

A virtual string is a counter-clockwise oriented copy of $S^1$, called a {\it core circle}, with arrows whose endpoints are glued to the core circle.  The endpoints of the arrows are required to be distinct.  We identify two virtual strings if there is a homeomorphism from one to the other preserving the directions of the arrows.

Every generic oriented curve on an oriented surface gives rise to a virtual string, called its {\it underlying virtual string}.  To construct the underlying virtual string of a curve, label the double points of the curve $a_1, \dots, a_n$.  Traverse the curve in the direction of its orientation and record the cyclic order in which the labels $a_i$ appear.  Each label appears twice in the cyclic order.  Then mark $2n$ points on a counterclockwise oriented copy of $S^1$, and label these points $a_1, \dots, a_n$ in the same cyclic order in which the labels appear on the double points during the traversal of the curve. For each $1\leq i \leq n$, connect the two points labelled $a_i$ by an arrow.  In other words, we connect the preimages of each double point of the map $S^1\rightarrow F$ by an arrow.  The direction of the arrows are determined by the following rule: At each intersection point of the curve order the two outgoing branches so that their tangent vectors form a positive frame.  The head of the arrow in the virtual string should point to the preimage corresponding to the first branch.  The underlying virtual string of a curve on a surface is also known as its Gauss diagram.

\begin{figure}[htbp]
\centering
\subfloat{
\includegraphics[width=3cm]{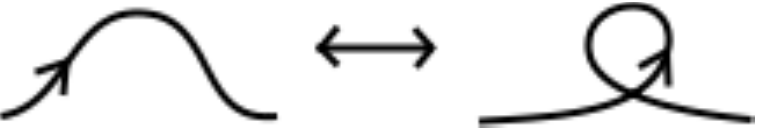}
\label{fig:subfig1}
}
\subfloat{
\includegraphics[width=3cm]{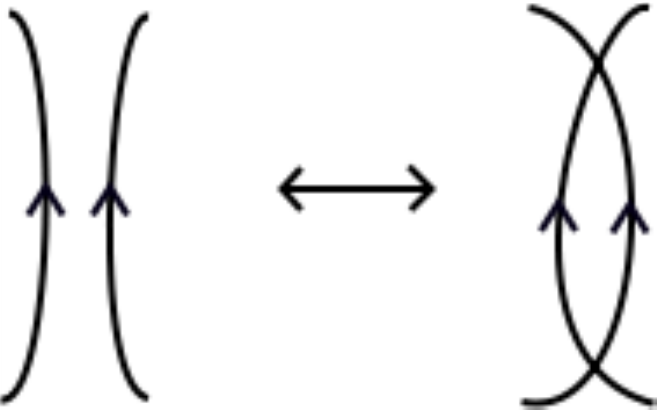}
\label{fig:subfig2}
}
\subfloat{
\includegraphics[width=3cm]{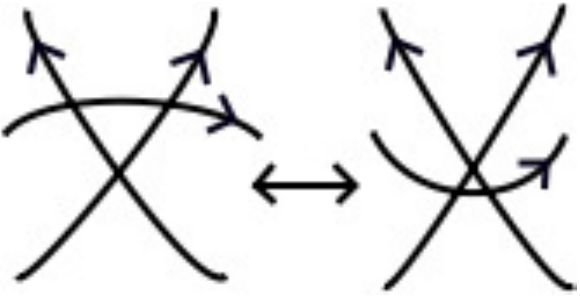}
\label{fig:subfig3}
}
\label{fig:subfigureExample}

\subfloat{
\includegraphics[width=3cm]{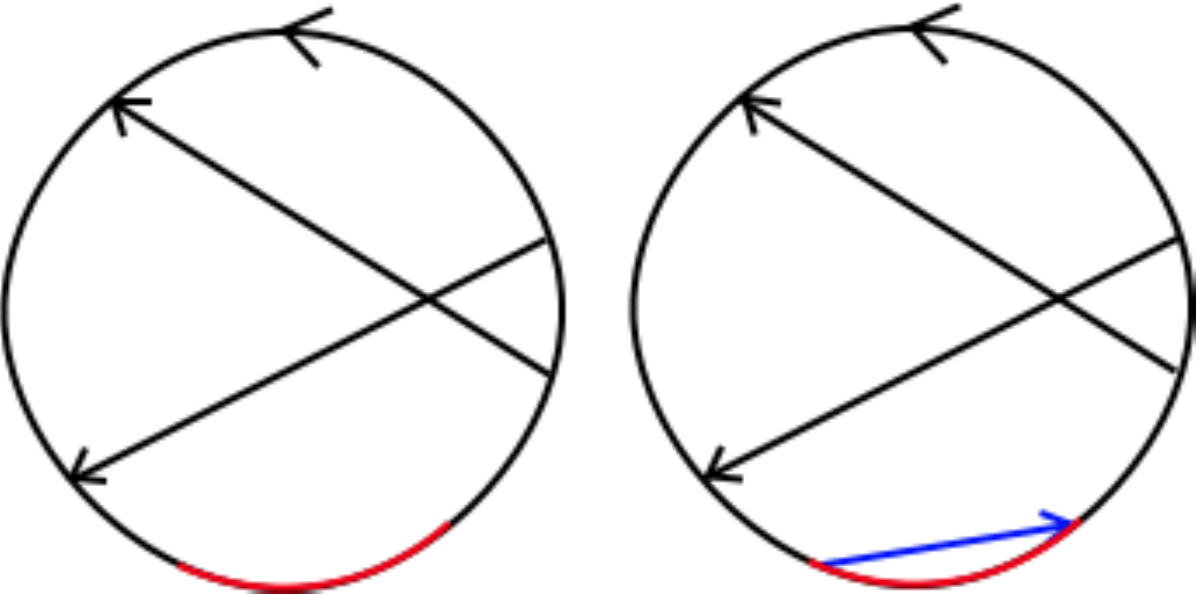}
\label{fig:subfig1}
}
\subfloat{
\includegraphics[width=3cm]{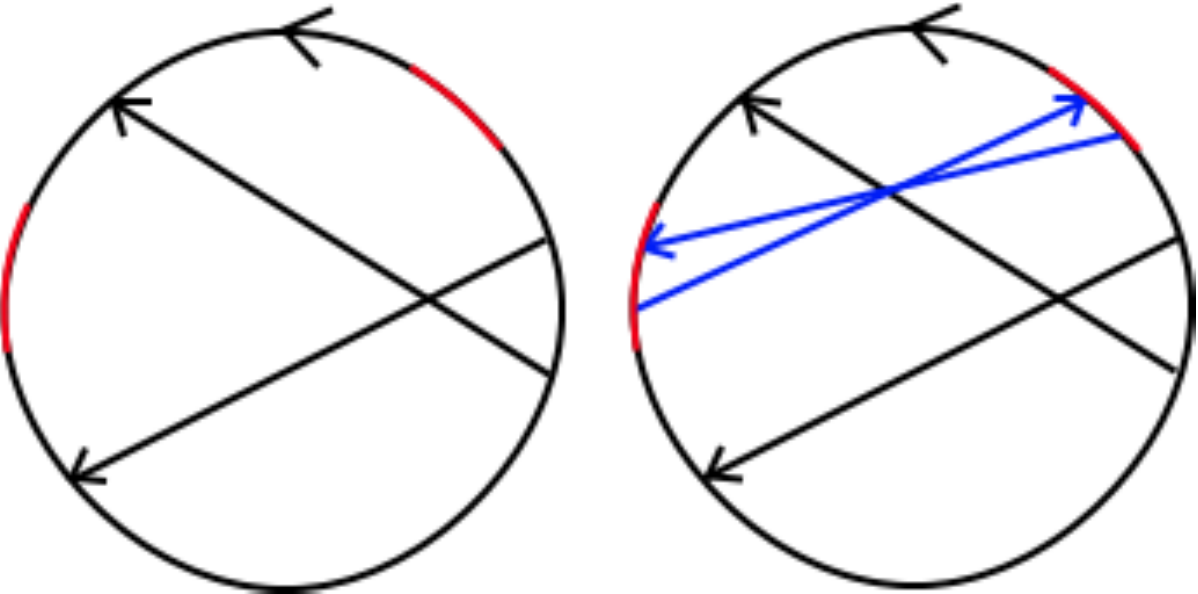}
\label{fig:subfig2}
}
\subfloat{
\includegraphics[width=3cm]{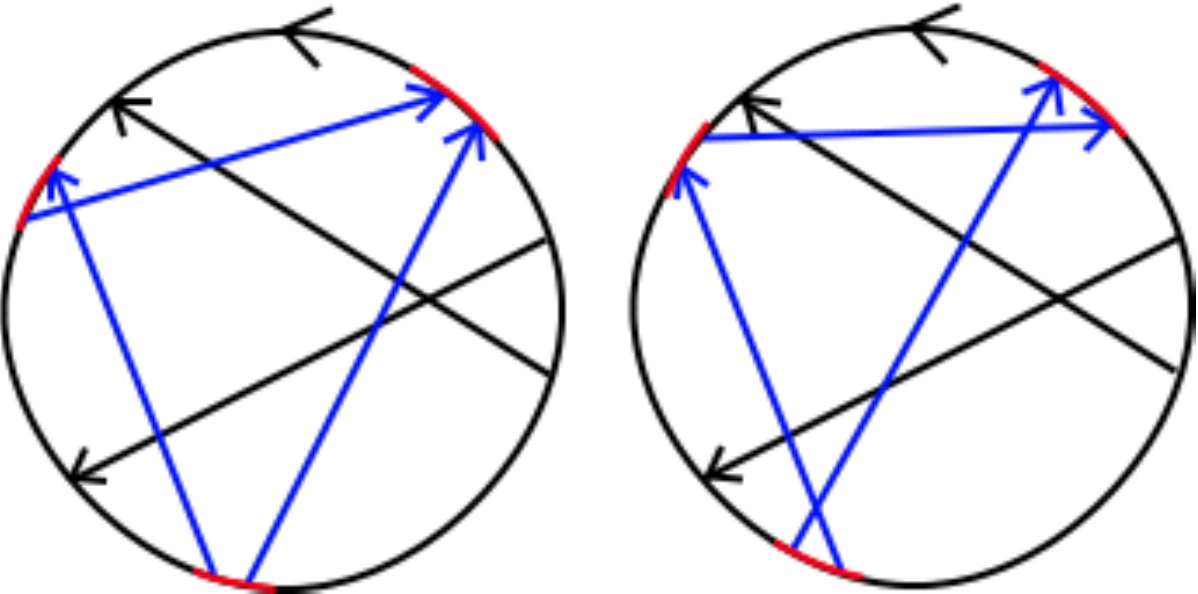}
\label{fig:subfig3}
}
\label{fig:subfigureExample}

\caption{Gauss diagram moves for flat virtual knots and their corresponding moves in the plane.}
\label{flatmoves.fig}
\end{figure}
Two curves on a surface are homotopic if and only if they are related by a sequence of flat Reidemeister moves, pictured at the top of Figure \ref{flatmoves.fig}, along with other moves that are obtained from those moves by independently reversing the orientation on any branch of the curve in the picture.  Each of these moves corresponds to a move on the underlying virtual string (Gauss diagram of the curve).  Two virtual strings are virtually homotopic if and only if they are related by a sequence of the Gauss diagram moves in the bottom row of Figure \ref{flatmoves.fig}, in addition to the Gauss diagram moves obtained from the other versions of the moves in the top row (which differ from the listed moves by the choice of orientation on the branches of the curve.)  In Section \ref{LegendrianGD} we will define Gauss diagram moves for virtual Legendrian knots.

\section{Gauss Diagrams for Virtual Legendrian Knots}\label{LegendrianGD}
We explain how to associate a Gauss diagram to an oriented virtual Legendrian knot.  For planar fronts, our diagrams are similar to the diagrams described by Polyak \cite{Polyak}.  However, unlike Polyak's diagrams, our diagrams are not marked with a basepoint and the signs on our cusps are different from Polyak's.

\begin{figure}[htbp]
	\includegraphics[width=8cm]{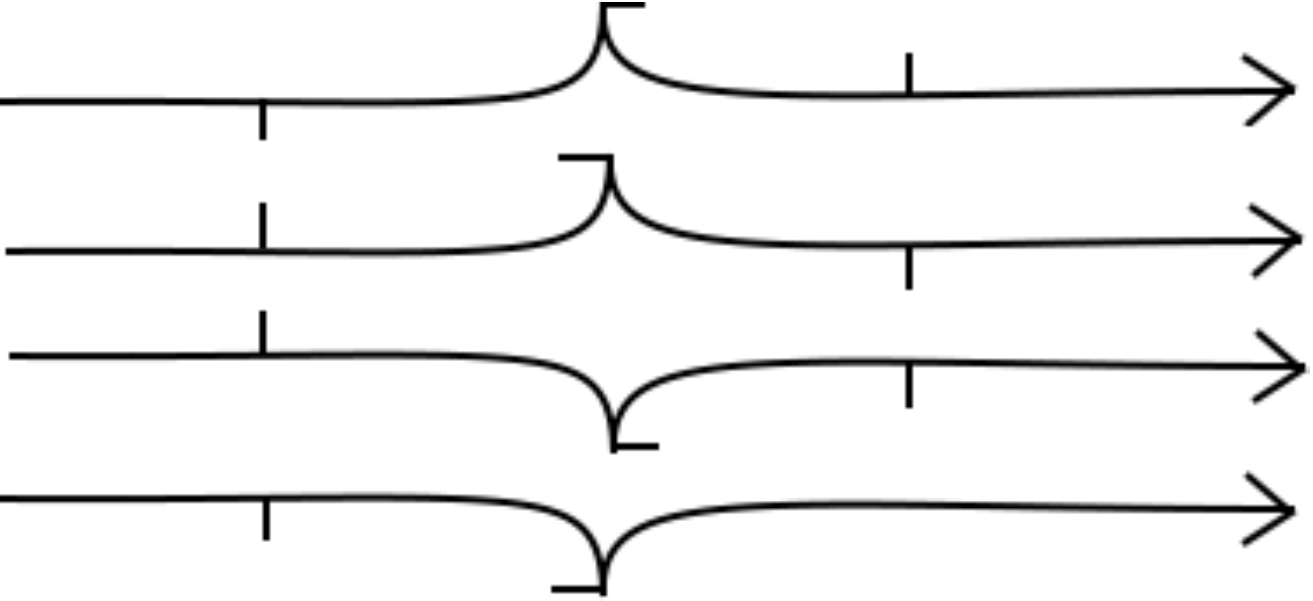}

	\caption{From top to bottom, a positive left cusp, negative left cusp, positive right cusp, and negative right cusp.}
	\label{cusps.fig}
\end{figure}

\begin{figure}[htbp]
	\includegraphics[width=8cm]{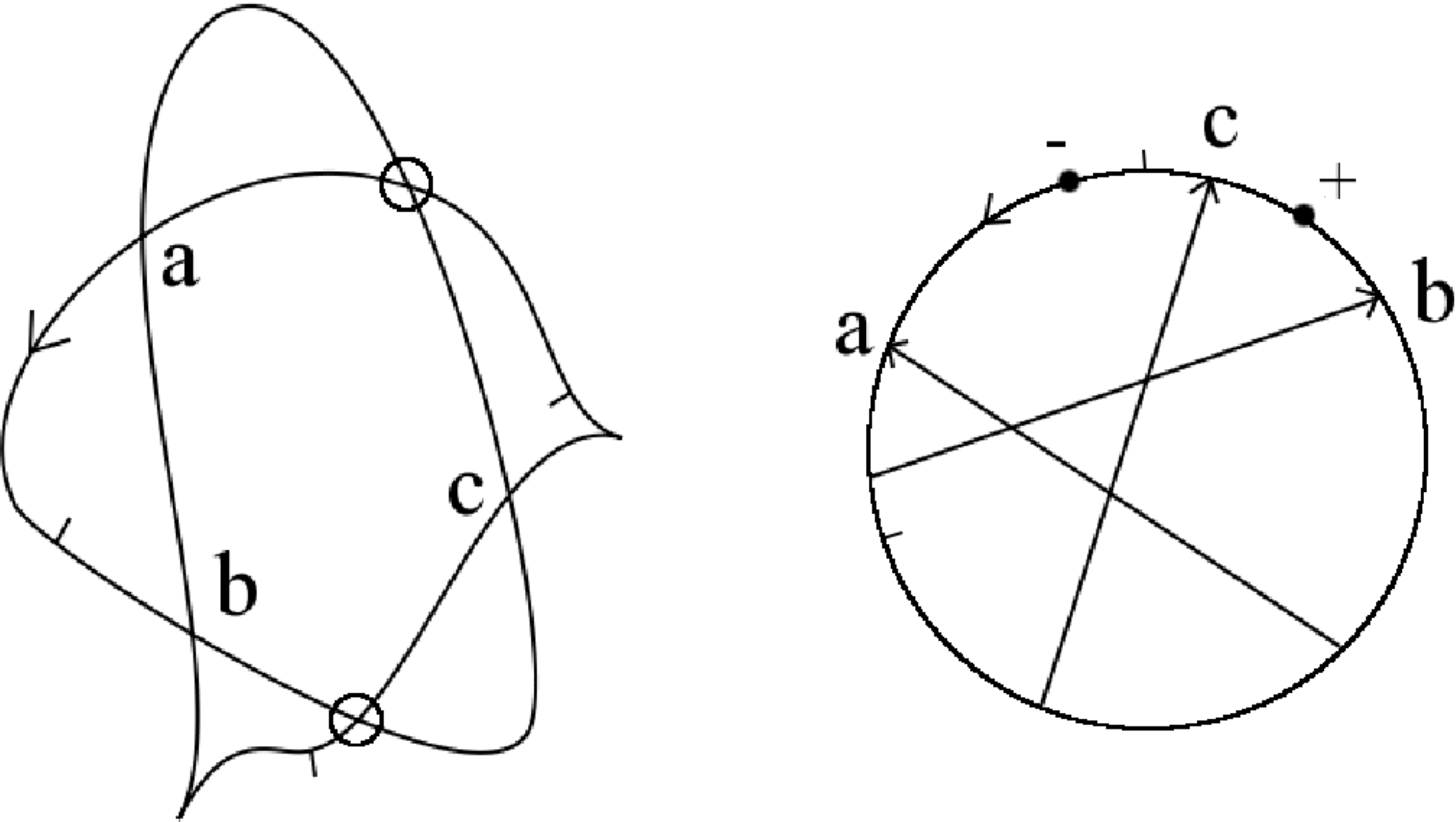}
	\caption{A virtual Legendrian knot diagram and its corresponding Legendrian Gauss diagram.}
	\label{Gaussdiagram.fig}
\end{figure}

The Gauss diagram of an oriented and cooriented wavefront $K$ on a surface $F$ with $c$ cusps and $n$ crossings (which we assume are transverse double points) is a counter-clockwise oriented copy of $S^1$, with $n$ arrows glued to $S^1$ at their endpoints, and $c$ marked points on $S^1$.  Let $C$ be the set of all marked points.  Each connected component of $S^1\setminus C$ is labeled with a coorientation.  Furthermore we require that the coorientations on adjacent components of $S^1\setminus C$ are different, and as a result $|C|$ is even.  The endpoints of the arrows are distinct, as are the marked points, and no marked point falls on the endpoint of an arrow.  The Gauss diagram of a given front is determined as follows.  We view $S^1$ as the circle parameterizing the curve $K$, and each pair of preimages of a double point of $K$ is connected by an arrow.  At each crossing, we label the outgoing branches of $K$ with a `1' and a `2' so that the ordered pair of their velocity vectors forms a positively oriented frame.  The head of the arrow is placed at the preimage corresponding to the branch labeled `1.'  The marked points of $C$ are the preimages of the cusps of $K$, and each marked point is equipped with a sign as follows.  A cusp of a wavefront is called positive (respectively, negative) if the outgoing branch of the cusp is (respectively, is not) in the coorienting half-plane of the cusp (see Figure \ref{cusps.fig}).  A virtual Legendrian knot diagram and its Gauss diagram are pictured in Figure \ref{Gaussdiagram.fig}.

Each wavefront move in Figure \ref{wavefrontmoves.fig} gives rise to a corresponding move on the Gauss diagram of the front in the obvious way, though we do not list these moves explicitly.  In Section \ref{equivalence} we will see that the equivalence relation generated by these Gauss diagram moves is simply the equivalence relation of virtual Legendrian knot theory.

\section{Virtual Framed and Virtual Legendrian Isotopy in the Spherical Cotangent Bundle}
The goal of this section is to define the notions of virtual framed and virtual Legendrian isotopy.  The definitions in this section are motivated by the definition of virtual knots in Carter, Kamada and Saito \cite{CKS} whose generalization to this context was suggested to us by Chernov.  We also introduce flat projection of a virtual blackboard framed knot. 

Throughout this section, $\bar{K}$ will typically denote a knot in $ST^*F$ and $K$ will typically denote its projection to the surface.

\subsection{The natural contact structure on the spherical cotangent bundle} \label{defSTF}
Let $F$ be an oriented surface and let $ST^*F$ be its spherical cotangent bundle.  That is, a point $\omega_p\in ST_p^*F$ is a linear functional on $T_pF$ defined up to multiplication by a positive scalar.  Hence $\omega_p$ is determined by a choice of 1-dimensional subspace $l_\omega$ of $T_pF$ such that $l_{\omega}=\ker \omega_p$, and a choice of positive half-space, which is a choice of connected component of $T_pF\setminus l_\omega$ on which $\omega_p$ is positive.  Put $\pi: ST^*F\rightarrow F$ to be the usual projection.  The contact plane at $\omega_p$ is $\pi_*^{-1}(l_\omega)$, which is a 2-dimensional subspace of $T_{\omega_p}ST^*F$.

\subsection{Virtual isotopy in the spherical cotangent bundle of a surface.}  Typically the virtual isotopy class of a knot is the isotopy class of a knot in a thickened surface $F\times I$ up to stabilization and destabilization of $F$.  In this paper we replace $F\times I$ with $ST^*F$.  We give a careful definition of virtual isotopy in $ST^*F$ in this section, and this definition is based on the formulation of virtual isotopy due to Carter, Kamada and Saito \cite{CKS}.

The {\it surface diagram} of a knot $\bar{K}:S^1\rightarrow ST^*F$ is a triple $(F,K,l)$, where $K$ is the projection of $\bar{K}$ to $F$, and $l$ is a cooriented line field along $K$ that describes how to lift $K$ to $\bar{K}$.  Namely, the point $K(t)$ lifts to the functional in $ST_{K(t)}^*F$ with kernel spanned by $l(t)$, and which is positive on the half-space of $T_{K(t)}F\setminus l(t)$ given by the coorientation of $l(t)$.

Now put $(F_1, K_1,l_1) \sim (F_2,K_2,l_2)$ if there exists a compact oriented surface $F_3$ and orientation preserving embeddings $\phi_1: F_1\rightarrow F_3$ and $\phi_2: F_2\rightarrow F_3$ such that the lifts of the surface diagrams $(F_3, \phi_1(K_1), \phi_{1*}(l_1))$ and $(F_3, \phi_2(K_2), \phi_{2*}(l_2))$ are isotopic in $ST^*F_3$.  Here $\phi_{i*}$ is the usual differential from $TF_i\rightarrow TF_3$, and $\phi_{i*}(l_i)$ is again a cooriented line field.  We abuse notation and also let $\phi_{i*}$ denote the natural map from $ST^*F_i$ to $ST^*F_3$.  At first glance it appears that this map goes in the incorrect direction, but $\phi_i$ is an embedding so its differential, and hence the induced map on cotangent spaces, are isomorphisms.  So, we abuse notation and let $\phi_{i*}:ST^*F_i\rightarrow ST^*F_3$. Note that the lift to $ST^*F_3$ of $(F_3,\phi_i(K_i),\phi_{i*}(l_i))$ is simply $\phi_{i*}(\bar{K_i})$.

Two virtual knots $(F,K,l)$ and $(F',K',l')$ are {\it virtually isotopic} if there  is a sequence of knot diagrams $(F_i,K_i,l_i)$, $1\leq i \leq m$, such that 
$$(F,K,l)=(F_1,K_1,l_1)\sim (F_2,K_2,l_2) \sim \dots \sim (F_m,K_m,l_m)=(F',K',l').$$

\subsection{Virtual framed isotopy in the spherical cotangent bundle} \label{framedisotopydef}
Next we define virtual isotopy for framed knots in $ST^*F$ where $F$ is an oriented surface.  A virtual framed knot is a knot $\bar{K}:S^1 \rightarrow ST^*F$ equipped with a transverse vector field $\nu$ considered up to the equivalence relation we define below.  We denote this virtual framed knot $(F,\bar{K}^\nu)$.  Because we will not need to work with the projection of a virtual framed knot, we will sometimes write $(F,K^\nu)$ rather than $(F,\bar{K}^\nu)$; the meaning of the notation will be clear from context.

Let $\phi: F\rightarrow F'$ be an orientation preserving embedding, and as above, let $\phi_*:ST^*F\rightarrow ST^*F'$ be the induced map on the spherical cotangent bundles. Let $\phi_{**}: TST^*F\rightarrow TST^*F'$ be the differential of $\phi_*$.

Put $(F_1, \bar{K_1}^{\nu_1}) \sim_f (F_2, \bar{K_2}^{\nu_2})$ if there exists a compact oriented surface $F_3$ and orientation preserving embeddings $\phi_1: F_1\rightarrow F_3$ and $\phi_2: F_2\rightarrow F_3$ such that $(F_3, \phi_{1*}(\bar{K_1})^{\phi_{1**}(\nu_1)})$ is framed isotopic to $(F_3, \phi_{2*}(\bar{K_2})^{\phi_{2**}(\nu_2)})$ in $ST^*F_3$.

  We say $(F_1, \bar{K_1}^{\nu_1})$ and $(F_m, \bar{K_m}^{\nu_m})$ are {\it virtually framed isotopic} if there exists a sequence of virtual framed knots $\bar{K_i}^{\nu_i}$, $1 < i < m$, satisfying
  $$(F_1, \bar{K_1}^{\nu_1})\sim_f (F_2, \bar{K_2}^{\nu_2}) \sim_f \dots \sim_f (F_m, \bar{K_m}^{\nu_m}).$$

\subsubsection{The blackboard framing}
We will sometimes consider knots with a certain framing, which we call the blackboard framing.  Fix an orientation of $ST^*F$ and of $F$.  A virtual topological knot $(F, K, l)$ is in general position if its velocity vector, $\bar{K}'(t)$, never points in the direction of the $S^1$ fiber.  The orientations of $ST^*F$ and $F$ determine an orientation of the $S^1$-fibers of $ST^*F$.  Let $\partial \theta$ be the vector field corresponding to the positive orientation of the $S^1$ fibers.  Then for a virtual topological knot in general position the vector field $\partial \theta$ is always transverse to $\bar{K}$ and thus is a framing vector field.  We call this framing the {\it blackboard framing}.  We can associate this framing to any virtual topological knot in general position to obtain a virtual framed knot with the blackboard framing, $\bar{K}^{\partial \theta}$.

\begin{figure}[htbp]
	\includegraphics[width=6cm]{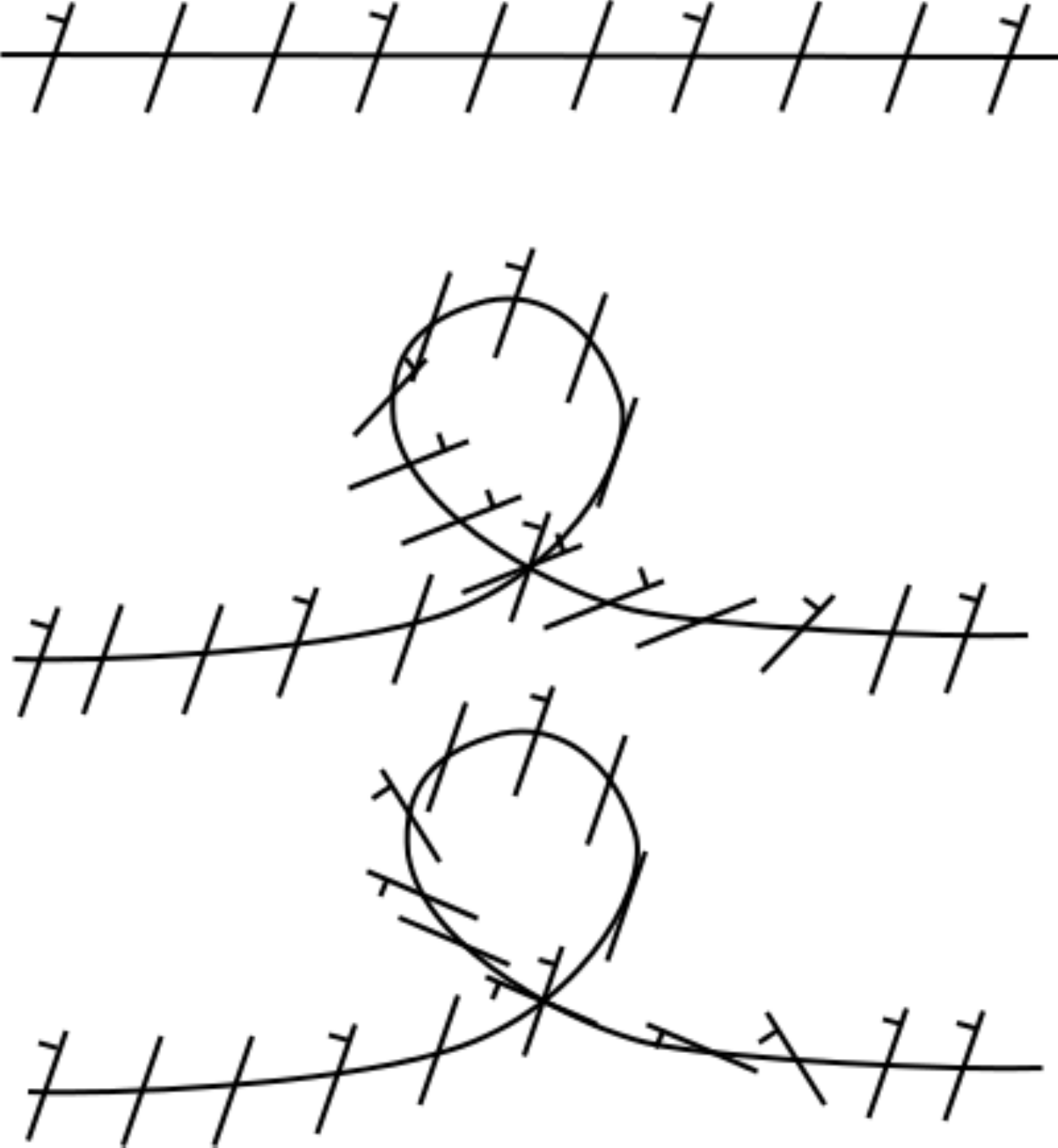}
	\caption{Adding a twist to the framing of a blackboard framed knot.}
	\label{blackboardkinks.fig}
\end{figure}

\begin{prop} Let $(F,\bar{K}^\nu)$ be a virtual framed knot.  Then there is a blackboard framed virtual knot $(F,\bar{L}^{\partial \theta})$ in the framed virtual isotopy class $[(F,\bar{K}^\nu)]_f$.
\end{prop}
\pp It is enough to show that the non-virtual framed isotopy class of $\bar{K}^\nu$ contains a blackboard framed knot $\bar{L}^{\partial \theta}$.  We may assume, possibly by first performing a small perturbation, that the velocity vector of $\bar{K}$ never points in the direction of the $S^1$ fiber, where $\bar{K}$ denotes the unframed knot in $ST^*F$ obtained from $\bar{K}^\nu$ by forgetting its framing.  Now consider the blackboard framed knot $\bar{K}^{\partial \theta}$ which coincides with $\bar{K}^\nu$ as an unframed knot.  Let $j=m(\bar{K}^\nu,\bar{K}^{\partial \theta})$ be the relative number of twists of the framings of the two knots (see Section ~\ref{mdefn}).  Next consider the surface diagram $(F,K,l)$ consisting of the projection $K$ of $\bar{K}$ to $F$ and the cooriented line field $l$ which describes how to lift $K$ to $\bar{K}$. We can replace a portion of this surface diagram, where the line field is locally constant, with a small kink (see Figure ~\ref{blackboardkinks.fig}).  Adding the top or bottom kink in Figure ~\ref{blackboardkinks.fig} to $(F,K,l)$ will yield new framed knots $\bar{K}_1$ and $\bar{K}_2$ respectively, which, once isotoped to coincide with $\bar{K}$ as embeddings, satisfy $m(\bar{K}_i,\bar{K})=\pm 1$, $i=1,2$.  The sign depends on the chosen orientation of the fiber.  We add $|j|$ copies of one of the kinks in Figure ~\ref{blackboardkinks.fig} so that the lift $\bar{L}$ of the resulting surface diagram is isotopic to $\bar{K}$, and the blackboard framed knot $\bar{L}^{\partial \theta}$ is framed isotopic to $\bar{K}^\nu$. 	\qed

\subsubsection{The flat projection of a virtual blackboard framed knot}
To prove Proposition $\ref{meven}$ we will use an invariant that is defined on the flat diagram of a blackboard framed virtual framed knot.  The flat projection of a virtual blackboard framed knot, $\bar{K}^{\partial \theta}$, is the immersed curve $\pi(\bar{K})$.

\begin{figure}[htbp]
	\includegraphics[width=8cm]{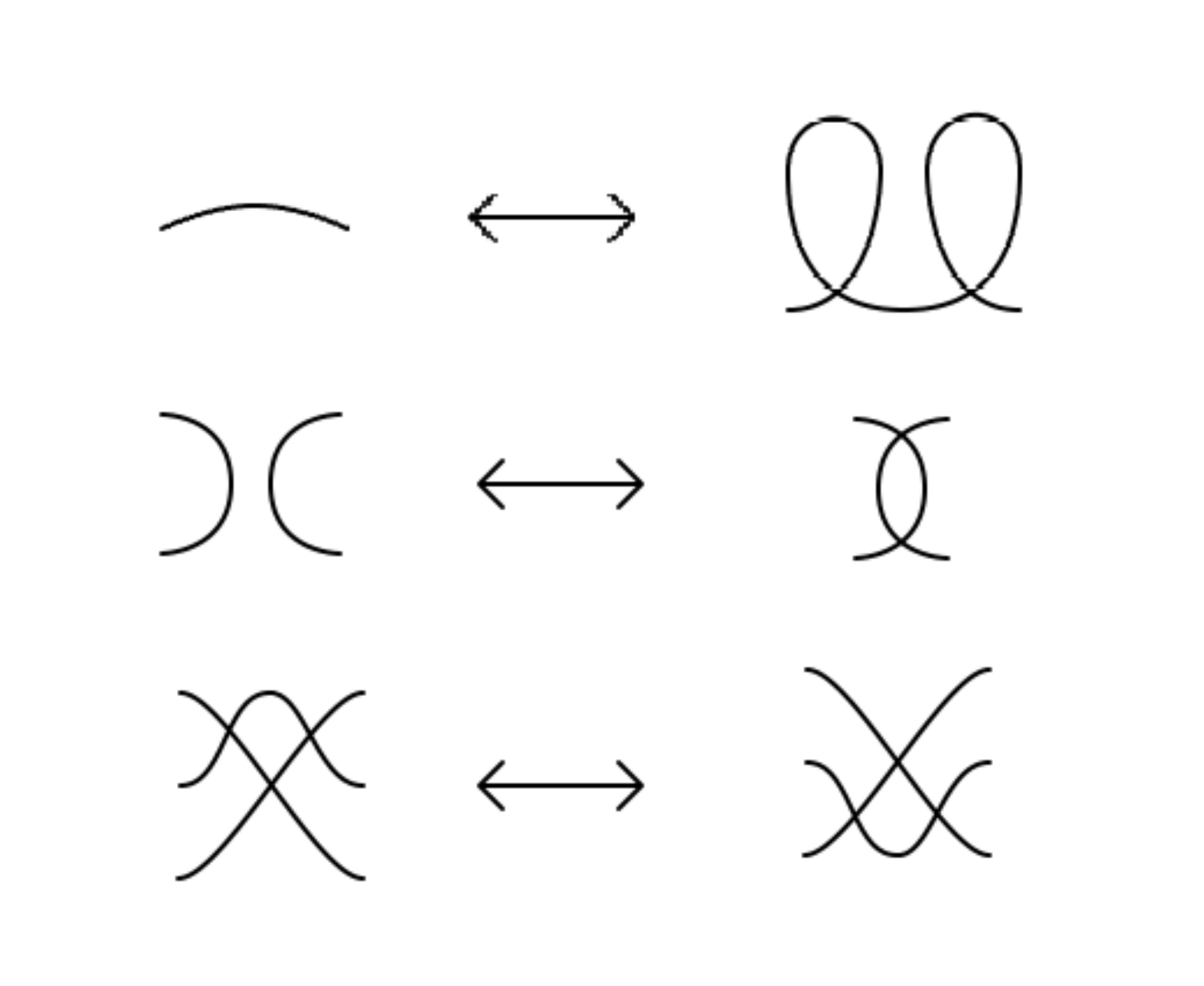}
	\caption{Moves for flat virtual blackboard framed knot projections}
	\label{flatBlackboardFrontMoves.fig}
\end{figure}

If two blackboard framed knots in the spherical cotangent bundle of the same surface are framed homotopic then their flat projections are related by a sequence of moves in Figure \ref{flatBlackboardFrontMoves.fig}.

Thus if two blackboard framed virtual knots $\bar{K}_1^{\partial \theta}$ and $\bar{K}_2^{\partial \theta}$ are virtually framed homotopic then there exists a sequence of moves in Figure \ref{flatBlackboardFrontMoves.fig}, in addition to stabilization and destabilization of the surface that takes the flat projection of $\bar{K}_1^{\partial \theta}$ to that of $\bar{K}_2^{\partial \theta}$.

\subsection{Virtual Legendrian Isotopy}\label{defvirtlegiso}
Let $\bar{K}$ be a Legendrian knot in $ST^*F$ in general position.  Then $\bar{K}$ projects to a cooriented wavefront $K$ on $F$.  Furthermore, two Legendrian knots $\bar{K}_1$ and $\bar{K}_2$ are Legendrian isotopic in $ST^*F$ if and only if their wavefronts $K_1$ and $K_2$ on $F$ are related by a sequence of the moves in Figure \ref{wavefrontmoves.fig}, excluding the dangerous self-tangency move.

Put $(F_1, K_1) \sim_l (F_2, K_2)$ if there exists a compact oriented surface $F_3$ and orientation preserving embeddings $\phi_1: F_1\rightarrow F_3$ and $\phi_2: F_2\rightarrow F_3$ such that $\phi_1(K_1)$ and $\phi_2(K_2)$ are related by a sequence of moves for wavefronts on $F_3$, or equivalently, if $\overline{\phi_1(K_1)}$ and $\overline{\phi_1(K_1)}$ are Legendrian isotopic in $ST^*F_3$.  

We say $(F,K)$ and $(F',K')$ are {\it virtually Legendrian isotopic} if there exists a sequence of pairs satisfying  $$(F,K)=(F_1,K_1)\sim_l (F_2, K_2)\sim_l \dots \sim_l (F_m,K_m)=(F',K').$$

\begin{rem} If in this definition we allowed dangerous tangency moves as well then we would have the definition of a virtual Legendrian homotopy.
\end{rem}

\begin{rem} \label{stdFraming} A virtual Legendrian knot in a cooriented contact structure has a natural framing, at each point given by the unit vector in the normal bundle to the velocity vector of the knot on which the coorienting one form evaluates to one.  Given a virtual Legendrian knot $\bar{K}$, let $\bar{K}^{st}$ denote the virtual framed knot given by $\bar{K}$ with this framing.
\end{rem}

\section{Flat Framed Virtual Knot Diagrams}\label{flatframed}
In this section we give reformulate the theory of flat virtual framed knots in terms of planar diagrams.

Given a virtual knot with a blackboard framing, $\bar{K}^{\partial \theta}$ we associate to it a planar flat virtual framed knot diagram.  This is obtained first by forgetting the framing and cooriented line field, leaving a generic immersed curve $K$ on the surface $F$, denoted $(F, K)$. Then we construct the  the virtual string associated to this curve as described in Section \ref{virtualstrings}.  This virtual string describes a flat virtual knot diagram in the plane in the usual way, which is unique up to any combination of virtual moves (a sequence of virtual moves is sometimes called the detour move).

\begin{figure}[htbp]
	\includegraphics[width=10cm]{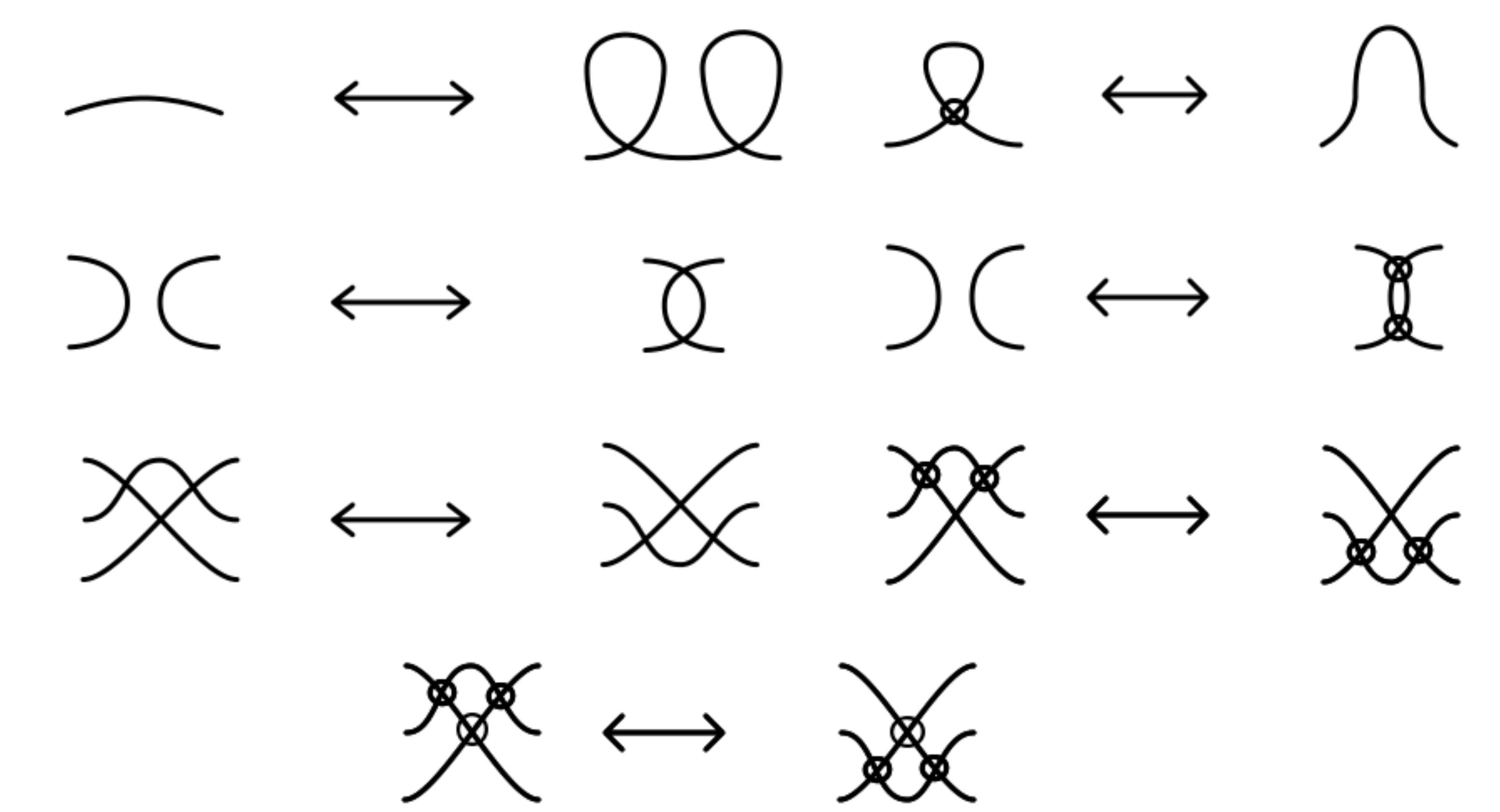}
	\caption{Moves for flat virtual framed knot diagrams}
	\label{flatDiagramMoves.fig}
\end{figure}

If two virtual framed knots are virtual framed homotopic then their planar flat virtual diagrams are related by a sequence of moves in Figure \ref{flatDiagramMoves.fig}.

\section{Equivalent Definitions of Virtual Legendrian Isotopy}\label{equivalence}
 Let $\mathcal{LGD}$ be the set of Legendrian Gauss diagrams and let $\sim_{lgd}$ be the equivalence relation generated by Gauss diagram moves.  Recall that the set of Gauss diagram moves is precisely the set of moves on Gauss diagrams obtained by translating each wavefront move to the Gauss diagram. Let $LGD=\mathcal{LGD}/ \sim_{lgd}$.  

Let $\mathcal{VLD}$ be the set of virtual Legendrian knot diagrams.  Let $\sim_{vld}$ be the equivalence relation generated by virtual Legendrian knot diagram moves.  Put $VLD=\mathcal{VLD}/\sim_{vld}$.

The theory of Legendrian Gauss diagrams up to Gauss diagram moves is equivalent to the theory of virtual Legendrian knot diagrams up to virtual wavefront moves.

\begin{thm}  The map $g: \mathcal{VLD}\rightarrow \mathcal{LGD}$ given by associating a (unique) Legendrian Gauss diagram to a virtual Legendrian knot diagram induces a bijection $g_\sim :VLD\rightarrow LGD$.
 \end{thm}
\pp  
First we verify that $g_\sim$ is well-defined.  Indeed, the Legendrian Gauss diagrams of two equivalent virtual Legendrian knot diagrams differ by Gauss diagram moves, as moves involving virtual crossings do not affect the Gauss diagram.  

The map $g_\sim$ is clearly surjective, as any Legendrian Gauss diagram gives rise to a virtual Legendrian knot diagram.  One simply draws all cusps and crossings present in the Legendrian Gauss diagram in the plane, and connects such cusps and crossings by arcs, creating virtual crossings where these arcs cross.

It remains to check that $g_\sim$ is injective.  Suppose we have two virtual Legendrian knot diagrams $D_1$ and $D_2$ with the same Gauss diagram.  We will show that $D_1$ and $D_2$ are related by virtual Legendrian knot diagram moves.  We first change $D_2$ by a regular isotopy so that a small neighborhood of each of its cusps and double points coincides with a small neighborhood of each corresponding cusp and double point of $D_1$.  The resulting diagrams differ only in how the cusps and double points are connected by arcs.  To move an arc $a_2$ of $D_2$ so that it coincides with the corresponding arc $a_1$ of $D_1$, we move $a_2$ by a fixed endpoint homotopy, such that any crossings created during that homotopy are virtual.  This move is known as the {\it detour move}, and is simply a sequence of the moves pictured in Figure \ref{virtualwavefrontmoves.fig}. 
\qed

Now let $\mathcal{SFD}$ be the set of all front diagrams on orientable surfaces, i.e., all pairs $(F,K)$ where $F$  is an oriented surface and $K$ is a cooriented wavefront on $F$.  Let $\sim_{sfd}$ be the equivalence relation generated by the relation $\sim_l$ defined in Subsection \ref{defvirtlegiso}, and let $SFD=\mathcal{SFD}/\sim_{sfd}$.  In other words, $SFD$ consists of pairs $(F,K)$ of oriented surfaces $F$ with cooriented wavefront diagrams $K$ on $F$ considered up to wavefront moves and stabilization and destabilization of the surface $F$.

Next we show that the theory of Legendrian Gauss diagrams up to Gauss diagram moves is equivalent to the theory of fronts on surfaces up to wavefront moves.
\begin{thm} The map $h: \mathcal{SFD}\rightarrow \mathcal{LGD}$ that assigns a Legendrian Gauss diagram to a wavefront on an oriented surface induces a bijection $h_\sim : SFD \rightarrow LGD$.
\end{thm}
\pp  First we check that $h_\sim$ is well-defined.  That is, suppose we have two pairs $(F,K)$ and $(F', K')$, such that for some sequence of pairs $\{(F_i,K_i)\}_{i=1}^n$, we have
$$(F,K)=(F_1,K_1)\sim_l (F_2,K_2)\sim_l \dots (F_n, K_n)=(F',K').$$
We need to check that $(F,K)$ and $(F',K')$ have equivalent Legendrian Gauss diagrams, but since stabilization and destabilization do not affect the Legendrian Gauss diagram of a wavefront on a surface, this is clear.

Next we verify that $h_\sim$ is surjective.  To do this, we construct a wavefront diagram on a surface given a Legendrian Gauss diagram.   First we build the disk-band surface realizing the underlying flat virtual knot of the wavefront.  For a detailed explanation of this procedure see \cite{Turaev}.  Then we insert positive and negative cusps according to the markings on the Legendrian Gauss diagram.

Finally we verify $h_\sim$ is injective.  To do this we must show that if two pairs $(F,K)$ and $(F',K')$ have equivalent Gauss diagrams, then $(F,K)$ and $(F',K')$ are equivalent.  Again, the argument is completely analogous to the case of virtual strings, which is carefully described in \cite{Turaev}. 
\qed

\section{Virtual Versions of the Classical Invariants}\label{classical.sec}

In this section we define the Maslov number for virtual Legendrian knots.  We do not discuss virtual analogues of the Bennequin number in this work. The virtual Maslov number $\mu(K_v)$ where $K_v$ is a planar virtual front diagram, is defined to be the number of positive cusps minus the number of negative cusps; see Figure \ref{cusps.fig}.  Clearly, if $K_v$ corresponds to the front $(F,K)$ on a surface then $\mu(K_v)$ is equal to the (non-virtual) Maslov number of the front $K$ on $F$.

A {\it positive (respectively negative) stabilization} of the virtual Legendrian knot $K$ is obtained by inserting a pair of positive (respectively negative) cusps at any point along $K$.  Let $K^{n_1,n_2}$ denote the virtual Legendrian knot obtained from $K$ by applying $n_1$ positive stabilizations and $n_2$ negative stabilizations.  

\begin{figure}[htbp]\includegraphics[width=9cm]{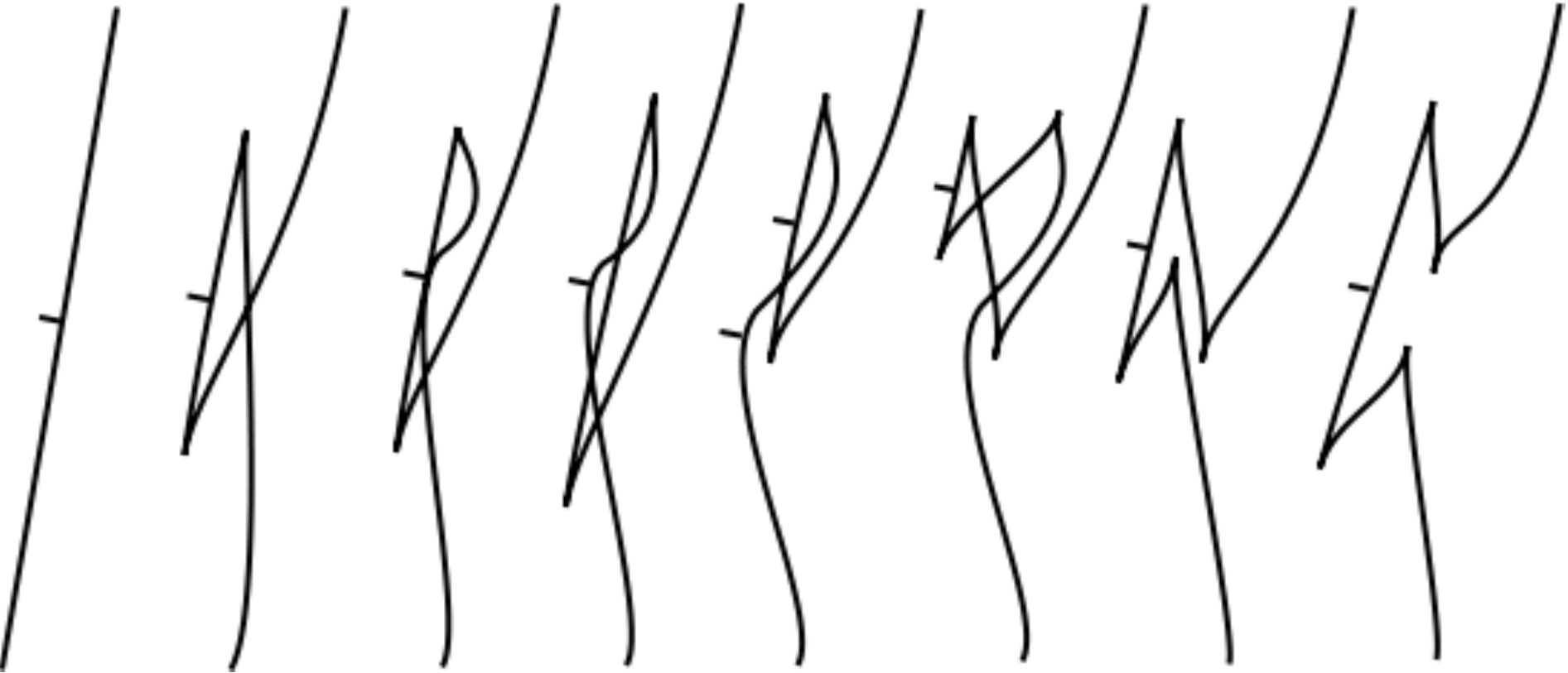}
	\caption{Adding a pair of negative cusps and a pair of positive cusps via a (virtual) Legendrian homotopy.}
	\label{homotopyposnegstabilization}
\end{figure}

\begin{prop}\label{propposnegstab} The virtual Legendrian knot $K$ is virtual Legendrian homotopic to $K^{n,n}$ for any positive integer $n$.
\end{prop}
\pp One can add a pair of positive cusps and a pair of negative cusps using a Legendrian homotopy; see Figure \ref{homotopyposnegstabilization}.  This sequence of moves is due to Fuchs and Tabachnikov \cite{f&t}.  This Legendrian homotopy is also a virtual Legendrian homotopy.
\qed

\begin{thm} Two virtual Legendrian knots are virtual Legendrian homotopic if and only if they have the same virtual Maslov number and are homotopic as virtual knots.
\end{thm}
\pp One can verify that the virtual Maslov number is an invariant of virtual Legendrian homotopy by checking its invariance under all virtual wavefront moves, as well as the dangerous tangency move.  Hence two virtual Legendrian homotopic virtual Legendrian knots have the same virtual Maslov number and (clearly) are homotopic as virtual knots.

Now suppose that $K$ and $L$ are Legendrian knots with the same virtual Maslov number that are homotopic as virtual knots.  We will show below that the assumption that $K$ and $L$ are homotopic as virtual knots implies that for any two sufficiently large positive integers $n_1$ and $n_2$, there exist integers $n_3$ and $n_4$ such that $K^{n_1, n_2}$ is virtual Legendrian homotopic to $L^{n_3,n_4}$.  In particular, this will be true for some $p=n_1=n_2$ large enough.  Then, since $K$ and $K^{p,p}$ are virtual Legendrian homotopic by Proposition \ref{propposnegstab}, and for suitable $n_3$ and $n_4$, $K^{p,p}$ and $L^{n_3,n_4}$ are virtual Legendrian homotopic, we have that $\mu(L)=\mu(K)=\mu(K^{p,p})=\mu(L^{n_3,n_4})$.  Then $\mu(L)=\mu(L^{n_3,n_4})$ implies $n_3=n_4$.  Put $q=n_3=n_4$.  Again by Proposition \ref{propposnegstab}, $L$ and $L^{q,q}$ are virtual Legendrian homotopic. Hence $K$ and $L$ are virtual Legendrian homotopic.

It remains to show that if $K$ and $L$ are virtually homotopic then for sufficiently large positive integers $n_1$ and $n_2$, there exist integers $n_3$ and $n_4$ such that $K^{n_1, n_2}$ is virtual Legendrian homotopic to $L^{n_3,n_4}$.  We do this in the next lemma. \qed

\begin{lem} \label{extrazigzaghomotopy} Let $(F,K)$ and $(F',L)$ be two virtually homotopic Legendrian knots.  Then for sufficiently large positive integers $n_1$ and $n_2$ there exist integers $n_3$ and $n_4$ such that $K^{n_1, n_2}$ is virtual Legendrian homotopic to $L^{n_3,n_4}$.
\end{lem}
\pp
We let $\bar{K}$ be a Legendrian knot in $ST^*F$, let $\bar{L}$ be a Legendrian knot in $ST^*F'$, and we have a sequence of pairs $$(F,K)=(F_1,K_1),(F_2,K_2),\dots,(F_n,K_n)=(F',L)$$ of curves $\bar{K_i}$ in $ST^*F_i$ such that the cooriented line fields $(F_i K_i, l_i)$ on $F_i$ and $(F_{i+1}, K_{i+1}, l_{i+1})$ on $F_{i+i}$ lifting to $\bar{K_i}$ and $\bar{K}_{i+1}$ respectively can be realized as homotopic cooriented line fields on a third surface $G_i$ (meaning their lifts are homotopic in $ST^*G_i$).  

On each surface $G_i$, this homotopy looks locally (within a Darboux chart) like a sequence of Reidemeister moves and crossing changes. We show how to imitate this virtual homotopy by a virtual Legendrian homotopy by replacing topological Reidemeister moves with Legendrian Reidemeister moves, and by replacing topological crossing changes with Legendrian crossing changes.
 
\begin{figure}[htbp]
\includegraphics[width=3cm]{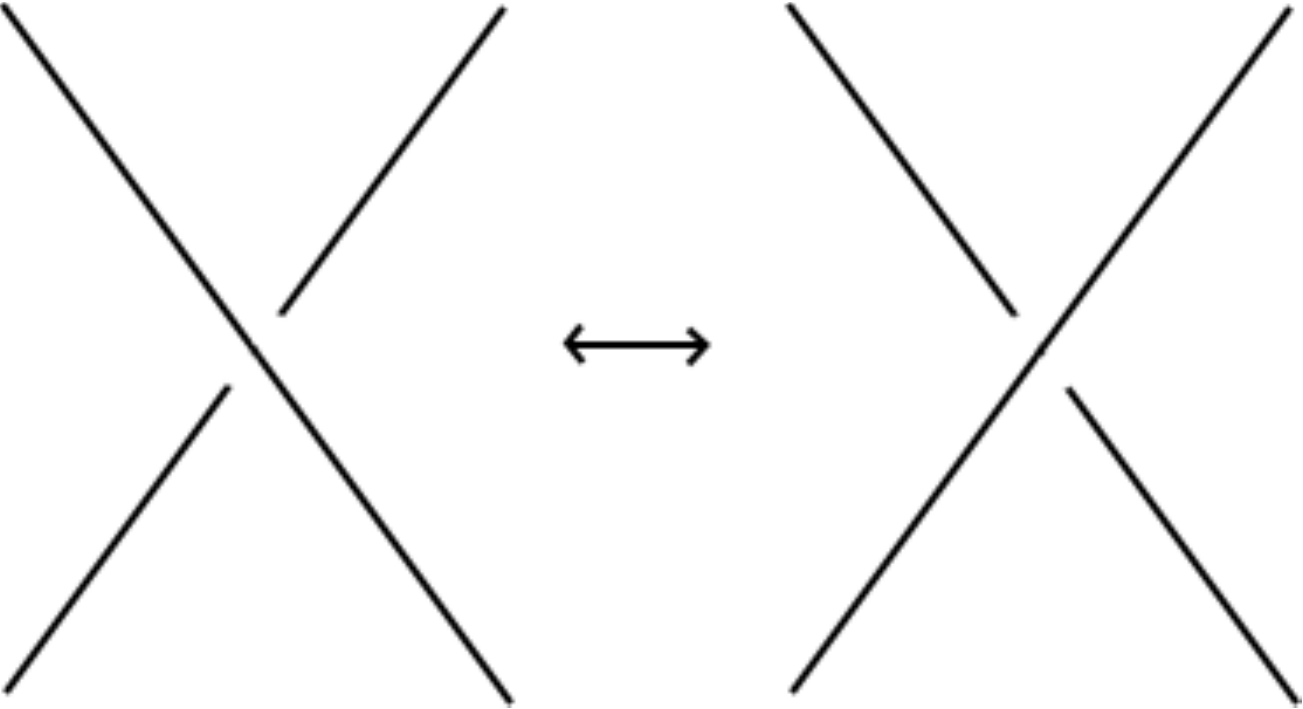}
\caption{A crossing change.}
\label{crossingchange.fig}
\end{figure}

The argument is now local, so we simply consider the case where $\bar{K}$ and $\bar{L}$ are homotopic as Legendrian knots in $ST^*F$ for a fixed surface $F$.  Furthermore we assume for now that the homotopy between $\bar{K}$ and $\bar{L}$ is contained in a single Darboux chart, so that $\bar{K}$ and $\bar{L}$ are Legendrian knots in the standard contact $\mathbb{R}^3$.  We consider the topological knot projections of $\bar{K}$ and $\bar{L}$ to the $xz$-plane.  Because $\bar{K}$ and $\bar{L}$ are homotopic, their topological projections are related by a sequence of Reidemeister moves of Types 1-3 and crossing changes; see Figure \ref{crossingchange.fig}.

Fuchs and Tabachnikov \cite{f&t} proved that if $K$ and $L$ are topologically isotopic Legendrian knots in the standard contact $\mathbb{R}^3$, then  for sufficiently large $n_1$ and $n_2$, there exist  $n_3$ and $n_4$ such that $K^{n_1, n_2}$ is Legendrian isotopic to $L^{n_3, n_4}$.  We prove the same statement, replacing isotopy with homotopy, and imitate their proof.  Let $\kappa$ and $\lambda$ be the front projections of $K$ and $L$. First we know $\kappa$ and $\lambda$ are related by a topological isotopy $K_t$ with projection $\kappa_t$, such that $K_0=K$, $K_1=L$, and for some $0=t_0<t_1< t_2 < \dots <t_n=1 \in [0,1]$, $\kappa_i=\kappa_{t_i}$ and $\kappa_{i+1}=\kappa_{t_{i+1}}$ are related by a single topological Reidemeister move, topological crossing change (see Figure \ref{crossingchange.fig}), or a passage through a vertical tangency (an ambient isotopy during which a strand without a vertical tangency passes through a vertical tangency, and after which two new vertical tangencies appear locally).

\begin{figure}[htbp]
\includegraphics[width=6cm]{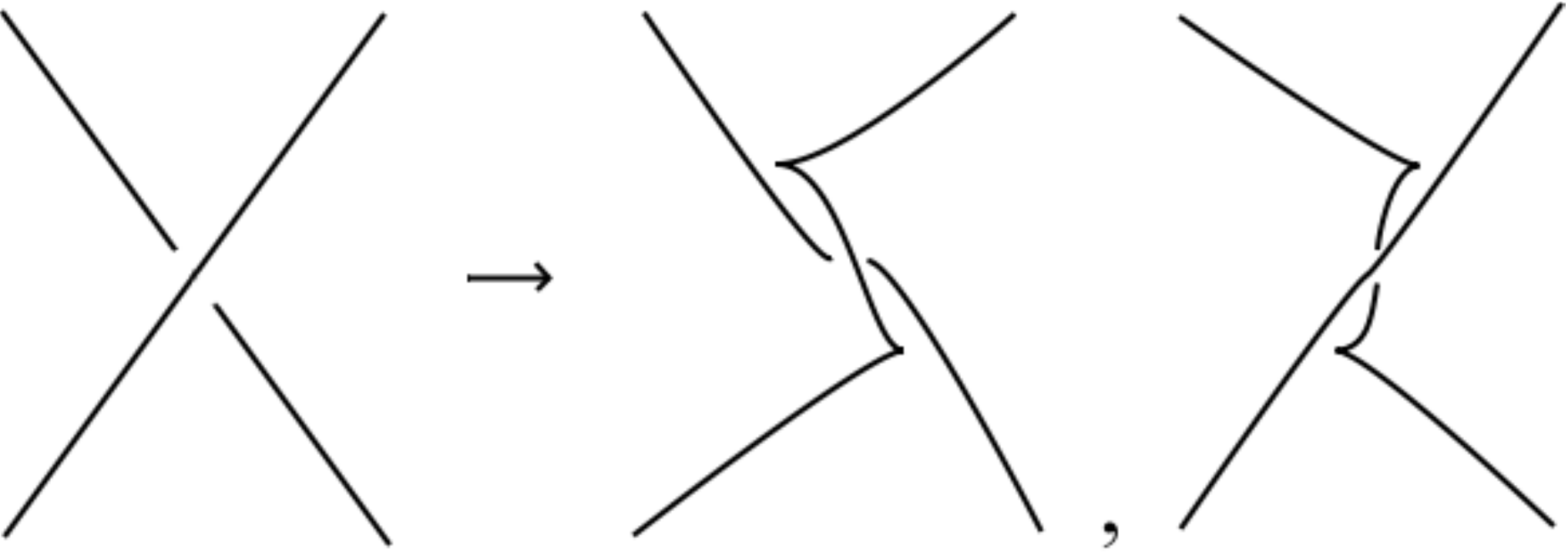}
\caption{Converting a ``wrong crossing'' in a topological knot diagram to a front projection of a Legendrian knot in standard contact $\mathbb{R}^3$.}
\label{wrongcrossing.fig}
\end{figure}

Next convert each topological knot diagram $\kappa_i$ into a front diagram by replacing vertical tangencies with cusps, and correcting ``wrong crossings" as in Fuchs and Tabachnikov \cite{f&t}, see Figure \ref{wrongcrossing.fig}.  These operations clearly do not affect the topological type of the knot.

\begin{figure}[htbp]
\includegraphics[width=9cm]{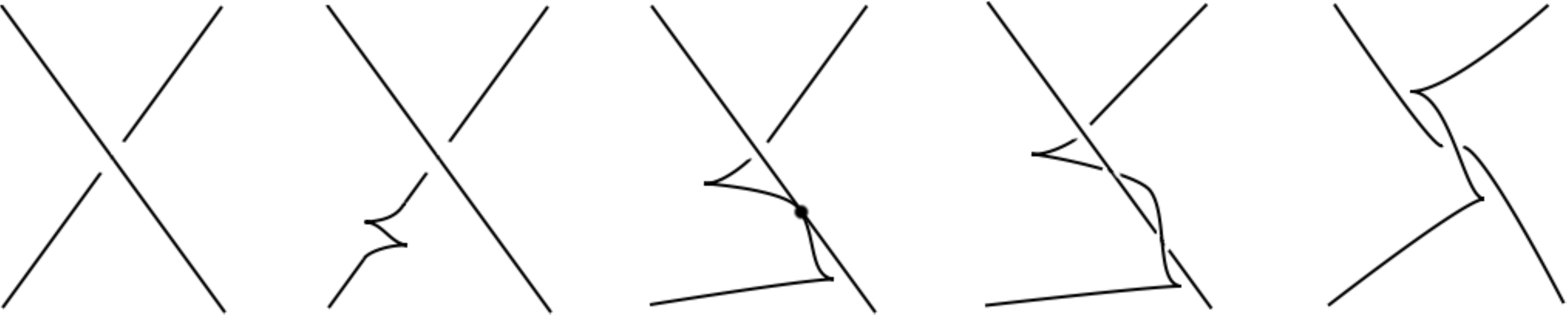}
\caption{After adding cusp pairs, we can replace a topological crossing change with a Legendrian crossing change.}
\label{legendriancrossingchange.fig}
\end{figure}

Fuchs and Tabachnikov then show how to replace topological Reidemeister moves and passages through vertical tangencies with Legendrian versions of these moves, possibly after adding extra positive or negative cusp pairs.  It remains to show that we can do the same for topological crossing changes.  This is shown in Figure \ref{legendriancrossingchange.fig}.
\qed

\section{Vassiliev Invariants}
In this section we define finite order Vassiliev invariants for virtual Legendrian, framed and topological knots. For the rest of the paper let $\mathcal{A}$ be an abelian group.

A self-intersection point of an immersed curve is called transverse if the two velocity vectors to the curve at that point are linearly independent.  An $n$-singular virtual Legendrian (resp. framed, topological) knot is a virtual Legendrian (resp. framed, topological) knot with $n$ transverse self-intersection points.

In an oriented $3$-manifold a transverse self-intersection point of an immersed curve can be resolved in two ways.  We call the resolution positive if the velocity vector to one strand, the velocity vector to the second strand and the vector from the second to the first strand form a positive $3$-frame.  Otherwise the self-intersection is called negative.  Given an immersed curve with $(n+1)$ transverse self-intersection points there are $2^{n+1}$ possible ways to resolve these self-intersection points. The {\it sign of a resolution} is the product of the signs of the resolutions of all of the individual double points. 

An $\mathcal{A}$-valued virtual Legendrian (resp. framed, topological) knot invariant is a function from the set of virtual Legendrian (resp. framed, topological) isotopy classes to $\mathcal{A}$.

A Vassiliev invariant of virtual Legendrian (resp. framed, topological) of order $\leq n$ is an $\mathcal{A}$-valued virtual Legendrian (resp. framed, topological) knot invariant that vanishes on the signed sum of the $2^{n+1}$ resolutions of any $(n+1)$-singular virtual Legendrian (resp. framed, toplogical) knot.
 
\section{Construction of the Isomorphism}
In what follows, recall that $[\cdot]_f$ denotes a framed virtual isotopy class, $[\cdot]_l$ denotes a Legendrian virtual isotopy class, and $[\cdot]$ denotes a topological virtual isotopy class.

When studying Vassiliev invariants of virtual Legendrian knots, we consider invariants of virtual Legendrian knots from a single connected component of the space of virtual Legendrian curves (or equivalently, knots within a fixed virtual homotopy class and with a fixed virtual Maslov number).   Similarly, when studying Vassiliev invariants of virtual framed knots or virtual Legendrian knots, we consider invariants of knots from a single connected component of the space of virtual framed curves.

In the following sections we will show that Vassiliev invariants cannot be used to distinguish virtually homotopic virtual Legendrian knots with the same virtual Maslov number that are isotopic as framed virtual knots.  First we prove the following seemingly weaker statement:

\begin{thm} \label{statementA}
Let $x\in \mathcal{V}_n^\mathcal{L}$.  Suppose $(F, K)$ and $(F', L)$ are two virtually homotopic virtual Legendrian knots with the same virtual Maslov number, such that $(F, K)$ is virtual framed isotopic to $(F', L)$ with their natural framings.   Then $x([(F, K)]_l) = x([(F', L)]_l)$.
\end{thm}

Using this theorem we are then able to prove the following stronger statement:

\begin{thm} \label{statementB}
Let $\mathcal{F}$ be a connected component of the space of virtual framed curves and $\mathcal{L}\subset\mathcal{F}$ be a connected component of the space of virtual Legendrian curves contained in $\mathcal{F}$.  Let $\mathcal{A}$ be an abelian group and $\mathcal{V}_n^\mathcal{F}$ be the group of $\mathcal{A}$ valued Vassiliev invariants on $\mathcal{F}$ of order $\leq n$.  Define $\mathcal{V}_n^\mathcal{L}$ likewise.  The restriction map $\phi:\mathcal{V}_n^{F}\rightarrow\mathcal{V}_n^\mathcal{L}$ is an isomorphism.
\end{thm}

It is clear that Theorem \ref{statementB} implies Theorem \ref{statementA}.  We will now show that Theorem \ref{statementA} implies Theorem \ref{statementB}.  We will prove Theorem \ref{statementA} in a later section.

Now we outline the construction of an inverse, $\psi$, to the map $\phi$ defined above.  If the framed isotopy class of $(F,K^\nu)$ contains a Legendrian representative $(F_l,K_l)$ then we define $\psi(x)([(F,K^\nu)]_f)=x([(F_l,K_l)]_l)$.  By Theorem \ref{statementA}, this is well-defined. If every virtual framed isotopy class were realizable by a virtual Legendrian knot, then the existence of $\psi$ would follow immediately from Theorem \ref{statementA}.  The Bennequin inequality tells us that not every (non-virtual) framed isotopy class is realizable by a Legendrian knot in, for example, $\mathbb{R}^3$ with the standard contact structure.  However, no Bennequin inequality is known for virtual Legendrian knots.  Therefore it is could be that every virtual framed isotopy class is realizable by a virtual Legendrian knot; this question is currently open.

\subsection{The relative number of twists of two framed knots.} \label{mdefn} Given two virtual framed knots $\bar{K}_1^{\nu_1}$ and $\bar{K}_2^{\nu_2}$ that coincide as pointwise embeddings in the same spherical cotangent bundle $ST^*F$, we can measure the relative number of twists of their framings as follows.  Let $\nu^{\perp}_1$ be a vector in $TST^*F$ orthogonal to both $\bar{K}'_1(t)$ and $\nu_1(t)$ such that the triple $\{\bar{K}'_1(t),\nu_1(t),\nu^{\perp}_1\}$ is a positive frame in $TST^*F$.  The frame consisting of $\nu^{\perp}_1$ and the vector ${\nu_1}(t)$ gives a trivialization of the normal bundle of $\bar{K}$.  Define $m(\bar{K}_1^{\nu_1},\bar{K}_2^{\nu_2})$ to be the total number of rotations of $\nu_2$ with respect to this trivialization.

In Proposition \ref{meven} we use $m$ to characterize when two framed virtual knots that coincide as unframed knots are homotopic as virtual framed curves.

For the proof of Proposition \ref{meven} we need the following definition.  Given a virtual framed knot $\bar{K}^\nu$ consider one of its corresponding flat virtual framed knot diagrams in $\mathbb{R}^2$ (see Section \ref{flatframed}).  Let $r$ be the rotation number of the flat framed knot diagram in the plane and let $v$ be the number of virtual crossings.  Put $\rho(\bar{K}^\nu) = r + v \mod{2}$.  One can check that this quantity does not change under any moves in Figure \ref{flatDiagramMoves.fig}, and thus is well-defined across all possible virtual framed knot diagrams for a given virtual framed knot. Furthermore this also shows that it is invariant under virtual framed homotopy. 

\begin{prop} \label{meven}
Let $\bar{K}_1^{\nu_1}$ and $\bar{K}_2^{\nu_2}$ be virtual framed knots (resp. singular virtual framed knots with n transverse double points) that coincide pointwise as embeddings (resp. immersons) in $ST^*F$.  Then $\bar{K}_1^{\nu_1}$ and $\bar{K}_2^{\nu_2})$ are virtual framed homotopic if and only if $m( \bar{K}_1^{\nu_1}, \bar{K}_2^{\nu_2})$ is even.
\end{prop}
\pp If $m( \bar{K}_1^{\nu_1}, \bar{K}_2^{\nu_2})$ is even, then $\bar{K}_1^{\nu_1}$ and $\bar{K}_2^{\nu_2}$ are homotopic as framed knots in $ST^*F$ because one can pass through a small kink to change $m$ by two.

Now suppose $m( \bar{K}_1^{\nu_1}, \bar{K}_2^{\nu_2})$ is odd.  Since $m$ is odd we must have that $\rho(\bar{K}_1^{\nu_1}) \neq \rho(\bar{K}_2^{\nu_2})$.  Thus $\bar{K}_1^{\nu_1}$ cannot be virtual framed homotopic to $\bar{K}_2^{\nu_2}$.
\hfill $\Box$

Suppose that $\bar{K}_1^{\nu_1}$ and $\bar{K}_2^{\nu_2}$ coincide as smooth embeddings, and $m(\bar{K}_1^{\nu_1}, \bar{K}_2^{\nu_2})=i$.  Then we write $\bar{K}_2^{\nu_2}=(\bar{K}_1^{\nu_1})^i$.

We want to prove that Theorem \ref{statementB} implies that the inverse $\psi$ described above exists. In \cite{Chernov} it is  shown that an analogue of Theorem \ref{statementB} for ordinary Legendrian and framed knots in most contact manifolds implies that an analogue Theorem \ref{statementA} holds, i.e., the inverse $\psi$ exists. The proof in \cite{Chernov} that $\psi$ exists is mostly local.  One can check that the same proof will work in the virtual category provided the following two propositions hold:

\begin{prop}  \label{prop2i}
Let $\mathcal{F}$ be a connected component in the space of virtual framed curves (resp. singular virtual framed curves), and $\mathcal{L}\subset\mathcal{F}$ be a connected component in the space of virtual Legendrian curves (resp. singular virtual Legendrian curves).  Let $\bar{K}^\nu\in\mathcal{F}$ be a virtual framed knot (resp. singular knot). Then there exists $i\in\mathbb{Z}$ and a virtual Legendrian knot (resp. singular knot) $\bar{L}\in \mathcal{L}$ such that $\bar{L}\in [(\bar{K}^\nu)^{2i}]_f$.  Furthermore if there exists a virtual Legendrian knot (resp. singular knot) $\bar{L}\in \mathcal{L}$ such that $[\bar{L}]_f = [\bar{K}^\nu]_f$ then there exists a virtual Legendrian knot (resp. singular knot) $\bar{L}'\in\mathcal{L}$ such that $[\bar{L}']_f = [(\bar{K}^\nu)^{-2}]_f$.
\end{prop}
\pp 

In \cite{Chernov} it is shown that for some $i \in \mathbb{Z}$, there exists a Legendrian knot $\bar{L}$ in the ordinary (non-virtual) framed isotopy class of $(\bar{K}^\nu)^{2i}$ that is also contained in the given connnected component of the space of Legendrian curves.  This $\bar{L}$ suffices here also because of Proposition ~\ref{meven}.  We will construct this knot $\bar{L}$ when the given knot is nonsingular but the singular case is similar and is done in \cite{Chernov} for non-virtual knots.  To construct $\bar{L}$ one first forgets the framing of the framed knot $\bar{K}^\nu$ to obtain a knot $K$. Then $C^0$-approximate the knot $K$ by a Legendrian knot, and add sufficiently many positive or negative cusp pairs, so that the resulting Legendrian knot $\bar{K}_l$ is in the given virtual Legendrian homotopy class $\mathcal{L}$.  This Legendrian knot $\bar{K}_l$ is framed isotopic to a framed knot $\bar{K}^{\nu j}$ such that $\bar{K}^{\nu j}$ and $\bar{K}^{\nu}$ coincide as unframed knots and $m(\bar{K}^{\nu j}, \bar{K}^{\nu})=j$.  But because $\mathcal{L}\subset \mathcal{F}$, then by Proposition ~\ref{meven}, $j=2i$.  

To show that if the framed isotopy class of $\bar{K}^\nu$ in $ST^*F$ is realizable by a Legendrian knot, then the framed isotopy class of $(\bar{K}^\nu)^{-2}$ in $ST^*F$ is realizable by a Legendrian knot $\bar{L}'$, simply preform the Legendrian homotopy in Figure ~\ref{homotopyposnegstabilization} to a small arc of the Legendrian knot $L$ in the virtual framed isotopy class of $\bar{K}^\nu$.  The resulting Legendrian knot is in the virtual framed isotopy class of $(\bar{K}^\nu)^{-2}$.
\qed

\begin{prop} \label{conncomp}  Let $\mathcal{F}$ be a connected component of the space of virtual framed curves, let $\bar{K}^\nu \in \mathcal{F}$ and let $K_u=(F,K,l)$ be an unframed virtual knot obtained by forgetting the framing on $\bar{K}^\nu$. Let $[K_u]$ be the class of virtual topological knots that contains $K_u$ and $\bar{K}_1^{\nu_1} = (F_1, \bar{K}_1^{\nu_1})\in \mathcal{F}$ be a virtual framed knot with $(F_1, K_1, l_1)\in [K_u]$.  Then $\bar{K}_1^{\nu_1}$ and $(\bar{K}^{\nu})^{2i}$ are virtual framed isotopic for some $i\in\mathbb{Z}$. 
\end{prop}
\pp This follows immediately from Proposition \ref{meven}.
\qed

We know how to define $\psi$ on virtual framed isotopy classes containing a virtual Legendrian knot. Let $\bar{K}^\nu \in \mathcal{F}$ and $i$ be the largest integer such that $[(\bar{K}^\nu)^{2i}]_f$ contains a virtual Legendrian knot in $\mathcal{L}$.  (If no such $i$ exists, then every $[(\bar{K}^\nu)^{2i}]_f$ contains a virtual Legendrian knot, so there is no problem defining $\psi$ on all of $\mathcal{F}$.) In this situation we know how to define $\psi$ on the framed isotopy classes $[K^{2j}]_f$ for $j\leq i$.  The following definition, analogous to the definition in \cite{Chernov}, extends $\psi$ to the virtual framed isotopy classes $[K^{2j}]_f$ for $j>i$. 

\begin{defin}
Fix $\bar{K}^\nu \in\mathcal{F}$ and let $j$ be the maximal integer such that $[(\bar{K}^\nu)^{2j}]$ contains a virtual Legendrian knot in $\mathcal{L}$.  For $l>j$ define 
$$\psi(x)((\bar{K}^\nu)^{2l}) = \sum_{i=1}^{n+1}\left( (-1)^{i+1} \frac{(n+1)!}{i!(n+1-i)!}\psi(x)((\bar{K}^\nu)^{2l-2i}) \right)$$
\end{defin}

This definition extends $x$ such that it is a Vassiliev invariant of virtual framed knots of order $\leq n$ and also $\phi \circ \psi = id_{\mathcal{V}_n^\mathcal{L}}$  and $\psi \circ \phi = id_{\mathcal{V}_n^\mathcal{F}}$ as we wanted.  The proof in our case is directly analagous to one given in \cite{Chernov}.

\section{Proof of Theorem \ref{statementA}}
Fix a connected component $\mathcal{F}$ of the space of virtual framed curves and a connected component $\mathcal{L}$ of the space of virtual Legendrian curves such that $\mathcal{L}\subset \mathcal{F}$.

In this section we will prove the following theorem:

\begin{thm} \label{samevalues}
Let $x\in \mathcal{V}_n^\mathcal{L}$.  Suppose $(F, K)$ and $(F', L)$ are two virtual Legendrian knots in $ \mathcal{L}$, such that $(F, K)$ is virtually framed isotopic to $(F', L)$.   Then $x([(F, K)]_l) = x([(F', L)]_l)$.
\end{thm}
Recall from Section \ref{classical.sec} that the virtual Legendrian knot $(F, K^{n, m})$ is obtained by adding $n$ positive and $m$ negative cusp pairs to the diagram of $(F, K)$.  

A crucial tool in the proof of Theorem \ref{samevalues} is the following lemma:

\begin{lem}\label{zigzag} Let $(F,K)$ and $(F', L)$ be two virtual Legendrian knots that are isotopic as virtual framed knots, and let $x\in\mathcal{V}^\mathcal{L}_n$. Suppose there exists $p\in \mathbb{Z}$ such that $(F,K^{p,p})$ and $(F',L^{p,p})$ are virtually Legendrian isotopic.  Then $x([(F,K)]_l)=x([(F',L)]_l)$.
\end{lem}

In order to use the previous lemma we must first show that the integer $p$ in Lemma \ref{zigzag} exists whenever $(F,K)$ and $(F',L)$ are in the same connected component of the space of virtual Legendrian curves and are isotopic as virtual framed knots.

To do this, we show that for $n_1,n_2$ large enough, there exists $n_3,n_4$ such that $(F,K^{n_1,n_2})$ and $(F',L^{n_3,n_4})$ are virtually Legendrian isotopic; this holds without the assumptions that $K$ and $L$ are virtually framed isotopic and are homotopic as virtual Legendrian cuves.  Then we show that we can assume that $n_1+n_2=n_3+n_4$ (provided $K$ and $L$ are virtually framed isotopic) and $n_1-n_2=n_3-n_4$ (provided $K$ and $L$ are homotopic as virtual Legendrian cuves).  It will follow that $n_1=n_2=n_3=n_4$.

\begin{lem}\label{fuchstab}
Let $(F, K)$ and $(F', L)$ be virtually isotopic Legendrian knots.  Then there exist $n_1, n_2, n_3$ and $n_4$ such that $(F,K^{n_1,n_2})$ is virtually Legendrian isotopic to $(F', L^{n_3, n_4})$.
\end{lem}
\pp
The proof is the same as the proof of Lemma \ref{extrazigzaghomotopy}, except that because we are now using isotopy rather than homotopy, we do not need to consider crossing changes.\qed

\begin{thm}\label{sameframediso}
Let $(F, K)$ and $(F', L)$ be two virtual Legendrian knots in the same virtual framed isotopy class $\mathcal{F}$. Then given large enough $n_1, n_2, n_3, n_4 \in \mathbb{Z}$ so that $(F, K^{n_1,n_2})$ and $(F',L^{n_3, n_4})$ are isotopic as virtual Legendrian knots, we have that $n_1+n_2 = n_3+n_4$.
\end{thm}

\pp
Throughout the proof below we write $K^\nu$ instead of $\bar{K}\nu$ to denote a framed knot in $ST^*F$ to increase readability.  $\bar{K}$ will still denote the lift to $ST^*F$ of a wavefront $K$.

Since $(F,K)$ and $(F',L)$ are virtual framed isotopic we have a sequence of pairs:
$$(F, K^{st}) = (F_1,K_1^{\nu_1}) \sim_f (F_2, K_2^{\nu_2}) \sim_f \dots \sim_f (F_m, K_m^{\nu_m}) = (F', L^{st}).$$

More precisely we have surfaces $G_i$ and maps $\phi_i : F_i \rightarrow G_i$ and $\psi_i:F_{i+1}\rightarrow G_i$ such that $\phi_{i*}(K_i^{\nu_i})$ is framed isotopic to $\psi_{i*}(K_{i+1}^{\nu_{i+1}})$ on $G_i$ via the framed isotopy $h_t^i$ for all $1\leq i < n$.

Now, by the argument in Lemma \ref{fuchstab} (or rather, Lemma \ref{extrazigzaghomotopy}) we can approximate the previous framed isotopy with a Legendrian isotopy after adding sufficiently many positive and negative cusp pairs to $K$ and $L$.  This yields a sequence
$$(F, K^{n_1,n_2}) = (F_1, L_1)\sim_l (F_2, L_2) \sim_l \dots \sim_l (F_m, L_m)=(F',L^{n_3,n_4})$$

with the same surfaces $G_i$ and maps $\phi_i$ and $\psi_i$ as above. However, now we have a sequence of Legendrian isotopies $l_t^i:S^1\rightarrow ST^*G_i$ such that the image of $l_0^i$ is equal to $\phi_{i*}(\bar{L_i})$,  and the image of $l_1^i $ is equal to $ \psi_{i*}(\bar{L}_{i+1})$.  Furthermore for all $t\in[0,1]$, the image of $l_t^i$ is contained in a small torus around the image of $h_t^i$.

We can use the fact that both of these images are contained in a small torus  at each time $t$ to show that $n_1+n_2 = n_3+n_4$.

Given two homotopic (virtual) framed knots $K_1^{\nu_1}$ and $K_2^{\nu_2}$ in $ST^*F$ which lie in a solid torus $T$, where $T$ is embedded in $ST^*F$, we can identify $T$ with the standard torus in $\mathbb{R}^3$.  Then we define $\text{slkd}(K_1^{\nu_1},K_2^{\nu_2})$ to be the difference of the self-linking numbers of the images of $K_1^{\nu_1}$ and $K_2^{\nu_2}$ under this identification.  Since both $K_1$ and $K_2$ are homotopic to the longitude of the torus, $\text{slkd}(K_1^{\nu_1}, K_2^{\nu_2})$ does not depend on the choice of identification.  Furthermore, one can check that $\text{slkd}(\bar{K},\bar{K}^{n_1,n_2})=n_1+n_2$, and $\text{slkd}(\bar{L},\bar{L}^{n_3,n_4})=n_3+n_4$. This argument is similar to an argument in \cite{Chernov}.

We now know that $\text{slkd}$ does not change as $t$ varies on a fixed surface $G_i$.  Thus we have that $\text{slkd}(\phi_{i*}(K_i^{\nu_i}), \phi_{i*}(\bar{L}_i)) = \text{slkd}(\psi_{i*}(K_{i+1}^{\nu_{i+1}}), \psi_{i*}(\bar{L}_{i+1}))$.  So to finish the proof we just need to show the following: 
$$\text{slkd}(\psi_{i*}(K_{i+1}^{\nu_{i+1}}), \psi_{i*}(\bar{L}_{i+1})) = \text{slkd}(\phi_{i+1*}(K_{i+1}^{\nu_{i+1}}), \phi_{i+1*}(\bar{L}_{i+1}))$$
However, since $\text{slkd}$ does not depend on the identification of the torus in $G_i$ or $G_{i+1}$ with the standard torus in $\mathbb{R}^3$ this equality is clear.
\qed

\begin{thm}\label{samemaslov}
Let $(F, K)$ and $(F', L)$ be two virtual Legendrian knots in the same virtual isotopy class that are homotopic as virtual Legendrian curves. If $n_1, n_2, n_3, n_4 \in \mathbb{Z}$ are large enough, so that $(F, K^{n_1,n_2})$ and $(F',L^{n_3, n_4})$ are isotopic as virtual Legendrian knots, then $n_1-n_2 = n_3-n_4$.
\end{thm}

\pp
From Proposition \ref{fuchstab} we have $n_1, n_2, n_3, n_4\in\mathbb{Z}$ so that $(F, K^{n_1,n_2})$ and $(F',L^{n_3, n_4})$ are isotopic as virtual Legendrian knots.  Since $(F, K)$ and $(F', L)$ are homotopic as virtual Legendrian curves we have that $\mu(F, K)=\mu(F', L)$.  For the same reason we get that $\mu(F, K^{n_1, n_2})=\mu(F', L^{n_3,n_4})$.  Finally since $(F, K^{n_1,n_2})$ is obtained from $(F, K)$ by adding $n_1$ upward cusp pairs and $n_2$ downward cusp pairs we have $\mu(F, K^{n_1,n_2})-\mu(F, K) = n_1 - n_2$.  Similarly we have $\mu(F', L^{n_3,n_4})-\mu(F', L)=n_3-n_4$.  So from the equalities above we can conclude that $n_1-n_2=n_3-n_4$
\qed

In order to prove the final theorem we need the following combinatorial lemma:

\begin{lem} \label{comboLem}
For $0\leq i < p$,
$$\sum_{k=\lceil\frac{i}{z+1}\rceil}^p (-1)^{k+1}\binom{p}{k}\binom{k(z+1)}{i} = 0$$
\end{lem}
\pp
We show this by comparing coefficients of a polynomial.
\begin{align*}
x^p\left(\sum_{k=1}^{z+1}(-1)^k\binom{z+1}{k}x^{k-1}\right)^p &= (1-(1-x)^{z+1})^p\\
 &= \sum_{k=0}^p (-1)^k\binom{p}{k}(1-x)^{k(z+1)}\\
 &= \sum_{k=0}^p (-1)^k\binom{p}{k}\left( \sum_{j=0}^{k(z+1)}(-1)^j\binom{k(z+1)}{j}x^j \right)\\
 &= \sum_{j=0}^{p(z+1)} \left( \sum_{k=\lceil\frac{i}{z+1}\rceil}^p (-1)^{k+j}\binom{p}{k}\binom{k(z+1)}{j} \right) x^j
\end{align*}
\qed

\begin{thm}\label{ppsameinvariant}
Let $x\in \mathcal{V}_n^\mathcal{L}$, and let $(F, K), (F', L) \in \mathcal{L}$.  Then if there exists $p$ such that  $(F, K^{p,p})$ and $(F', L^{p,p})$ are virtual Legendrian isotopic then $x([(F, K)]_l) = x([(F', L)]_l)$.
\end{thm}

\begin{figure}[htbp]
	\includegraphics[width=6cm]{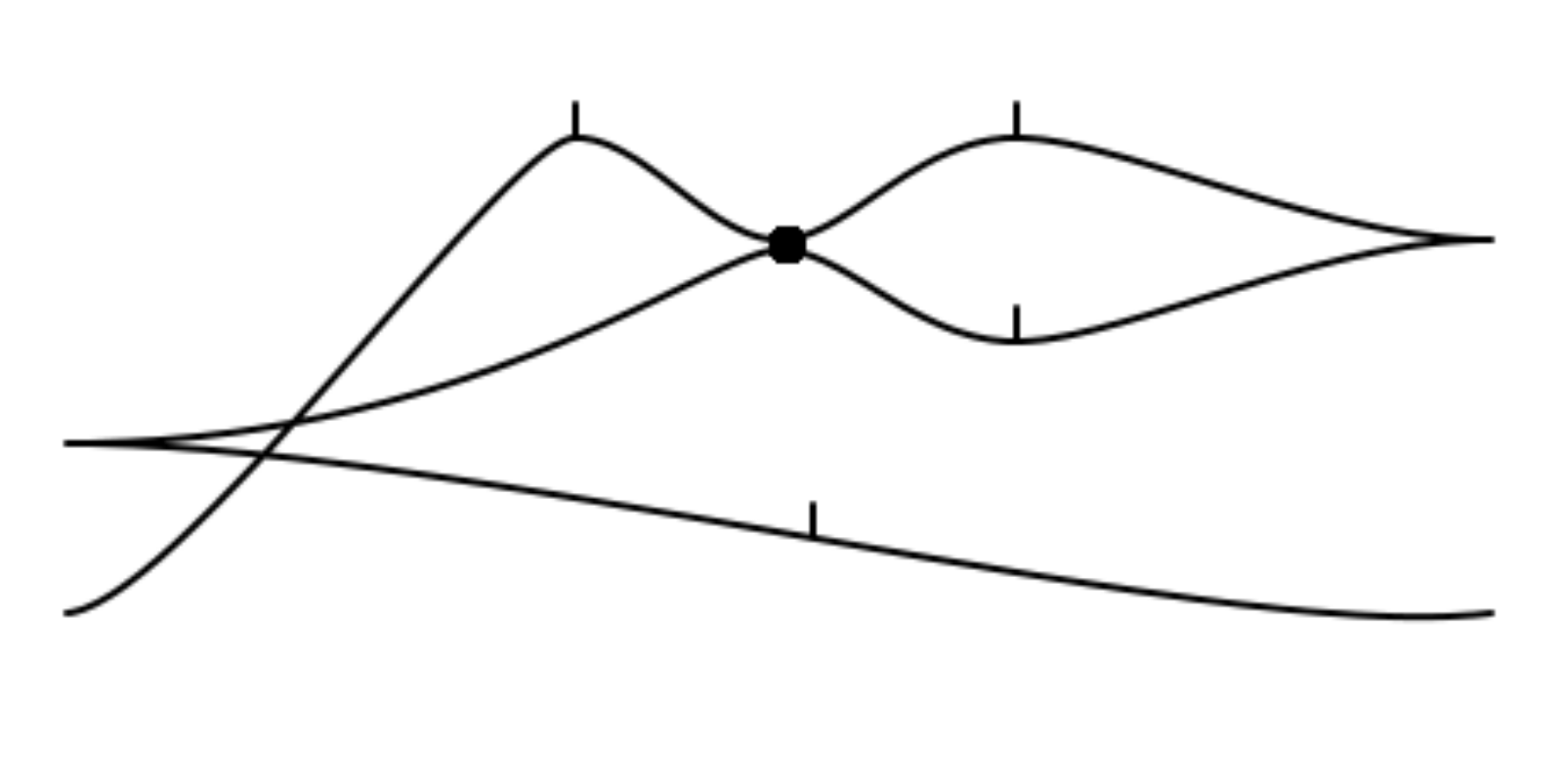}
	\caption{The singular Legendrian arc used to construct $K^z$.}
	\label{DoublePoint.fig}
\end{figure}

\pp Throughout the proof below we drop the $[\cdot]_l$ notation, with the understanding that $x$ is always an invariant of virtual Legendrian knots. Fix a point $q$ in the image of $(F, K)$ and denote by $(F, K^z)$ the singular virtual Legendrian knot obtained from $(F, K)$ by adding $z$ copies of Figure \ref{DoublePoint.fig} in a neighborhood of $q$.  For a singular knot $(F, K_s)$ denote by $d(F,  K_s)$ the sum of all the signed resolutions of the double points of $\bar{K}_s$, i.e., direct self-tangencies of $K_s$.  The homotopy in Figure ~\ref{homotopyposnegstabilization} implies that $d(F, K^1) = (F, K) - (F, K^{1,1})$.  So we have that $x(F, K) = x(F, K) - x(F, K^{1,1})$.  By iterating this process we can conclude that $x(F, K^z) = \sum_{j=0}^{z}(-1)^j\binom{z}{j}x\left(F, K^{j,j}\right)$.

To obtain the result we will begin with $x(F,K)$ and, after adding a linear combination of singular knots on which $x$ is zero, we will change the argument of $x$ to $(F', L)$.  First we make some observations that will be used in the computation:
\begin{enumerate} 
\item $x(F, K) = x(\sum_{k=0}^p(-1)^k\binom{p}{k}(F, K^{k(z+1)}))$, since for $k>0$, $(F, K^{k(z+1)})$ has at least $z+1$ double points.  
\item $x(F, K^{m,m}) = x(F', L^{m,m})$ for $m\geq p$ as you can change $(F, K^{p,p})$ to $(F', L^{p,p})$ through a virtual Legendrian isotopy
\item $x(F, K^z) = \sum_{j=0}^{z}(-1)^j\binom{z}{j}x\left(F, K^{j,j}\right)$ as noted above.
\end{enumerate}

\begin{align*}
x(F, K) &= \sum_{k=0}^p(-1)^k\binom{p}{k} x\left(F, K^{k(z+1)}\right)& \text{by (1) above}\\
&= \sum_{k=0}^p(-1)^k\binom{p}{k}\left(\sum_{j=0}^{k(z+1)}(-1)^j\binom{k(z+1)}{j}x\left(F, K^{j,j}\right)\right)& \text{by (3) above}\\
&= \sum_{j=0}^{p(z+1)} \left( \sum_{k=\lceil\frac{i}{z+1}\rceil}^p (-1)^{k+j}\binom{p}{k}\binom{k(z+1)}{j} \right) x\left(F, K^{j,j}\right)&\\
&= \sum_{j=p}^{p(z+1)} \left( \sum_{k=\lceil\frac{i}{z+1}\rceil}^p (-1)^{k+j}\binom{p}{k}\binom{k(z+1)}{j} \right) x\left(F, K^{j,j}\right)& \text{by } \ref{comboLem}\\
&= \sum_{j=p}^{p(z+1)} \left( \sum_{k=\lceil\frac{i}{z+1}\rceil}^p (-1)^{k+j}\binom{p}{k}\binom{k(z+1)}{j} \right) x\left(F', L^{j,j}\right)& \text{by  (2) above} \\
&= \sum_{j=0}^{p(z+1)} \left( \sum_{k=\lceil\frac{i}{z+1}\rceil}^p (-1)^{k+j}\binom{p}{k}\binom{k(z+1)}{j} \right) x\left(F', L^{j,j}\right)&\\
&= \sum_{k=0}^p(-1)^k\binom{p}{k}\left(\sum_{j=0}^{k(z+1)}(-1)^j\binom{k(z+1)}{j}x\left(F', L^{j,j}\right)\right)&\\
&= \sum_{k=0}^p(-1)^k\binom{p}{k}x\left(F', K_2^{k(z+1)}\right)&\text{by (3) above}\\
&= x(F', L)
\end{align*} 
\qed

Now by Theorems \ref{sameframediso}, \ref{samemaslov} and \ref{ppsameinvariant} we can conclude that Vassiliev invariants cannot distinguish virtually framed isotopic virtual Legendrian knots with the same Maslov number.
\begin{thm}
Let $x\in \mathcal{V}_n^\mathcal{L}$, $(F, K), (F', L) \in \mathcal{L}$ and let $(F, K)$ be virtually framed isotopic to $(F', L)$ then $x([F, K)]_l) = x([(F', L)]_l)$.
\end{thm}

{\it Acknowlegdement:} We would like to thank our advisor, Vladimir Chernov, for suggesting the problem and for many valuable discussions.

\end{document}